\input amssym.def
\input amssym.tex


\def\item#1{\vskip1.3pt\hang\textindent {\rm #1}}


\tolerance=300
\pretolerance=200
\hfuzz=1pt
\vfuzz=1pt


\hoffset=0.6in
\voffset=0.8in

\hsize=5.8 true in 
\vsize=8.5 true in
\parindent=25pt
\mathsurround=1pt
\parskip=1pt plus .25pt minus .25pt
\normallineskiplimit=.99pt

\countdef\revised=100
\mathchardef\emptyset="001F 
\chardef\ss="19
\def\3{\ss}
\def\anf{$\lower1.2ex\hbox{"}$}
\def\frac#1#2{{#1 \over #2}}
\def\>{>\!\!>}
\def\<{<\!\!<}

\def\into{\hookrightarrow}
\def\ssarr{\hbox to 30pt{\rightarrowfill}}
\def\sarr{\hbox to 40pt{\rightarrowfill}}
\def\arr{\hbox to 60pt{\rightarrowfill}}
\def\larr{\hbox to 60pt{\leftarrowfill}}
\def\Arr{\hbox to 80pt{\rightarrowfill}}

{}

\def\ad{\mathop{\rm ad}\nolimits}

\def\Ad{\mathop{\rm Ad}\nolimits}

\def\Aut{\mathop{\rm Aut}\nolimits}

\def\det{\mathop{\rm det}\nolimits}

\def\Diff{\mathop{\rm Diff}\nolimits}

\def\End{\mathop{\rm End}\nolimits}
\def\Ext{\mathop{\rm Ext}\nolimits}

\def\Gl{\mathop{\rm Gl}\nolimits}

\def\Herm{\mathop{\rm Herm}\nolimits}

\def\id{\mathop{\rm id}\nolimits} 
\def\im{\mathop{\rm im}\nolimits}
\def\Im{\mathop{\rm Im}\nolimits}

\def\Int{\mathop{\rm int}\nolimits}

\def\rank{\mathop{\rm rank}\nolimits}

\def\Re{\mathop{\rm Re}\nolimits}

\def\Sl{\mathop{\rm Sl}\nolimits}
\def\SO{\mathop{\rm SO}\nolimits}
\def\span{\mathop{\rm span}\nolimits}

\def\Symm{\mathop{\rm Symm}\nolimits}
\def\Sp{\mathop{\rm Sp}\nolimits}
\def\Spec{\mathop{\rm Spec}\nolimits}

\def\SU{\mathop{\rm SU}\nolimits}
\def\sup{\mathop{\rm sup}\nolimits}

\def\tr{\mathop{\rm tr}\nolimits}

\def\0{{\bf 0}}
\def\1{{\bf 1}}

\def\a{{\frak a}}

\def\e{{\frak e}}
\def\f{{\frak f}}
\def\g{{\frak g}}

\def\k{{\frak k}}
\def\l{{\frak l}}
\def\m{{\frak m}}

\def\n{{\frak n}}

\def\p{{\frak p}}
\def\q{{\frak q}}

\def\s{{\frak s}}

\def\sp{{\frak {sp}}}

\def\su{{\frak {su}}}
\def\so{{\frak {so}}}
\def\sL{{\frak {sl}}}
\def\t{{\frak t}}
\def\uu{{\frak u}}

\def\z{{\frak z}}

\def\C{{\Bbb C}} 
\def\D{{\Bbb D}} 
 
\def\H{{\Bbb H}}

\def\N{{\Bbb N}} 
 
\def\P{{\Bbb P}} 
 
\def\R{{\Bbb R}} 
 
\def\Z{{\Bbb Z}} 

\def\:{\colon}  
\def\.{{\cdot}}
\def\|{\Vert}
\def\bsk{\bigskip}

\def\giantskip{\vskip2\bigskipamount}
\def\gsk{\giantskip}
\def \la {\langle}
\def\msk{\medskip}
\def \ra {\rangle}
\def \res {\!\mid\!\!}

\def\ssk{\smallskip}

\def\bbr{\bigbreak}
\def\giantbreak{\par \ifdim\lastskip<2\bigskipamount \removelastskip
         \penalty-400 \giantskip\fi}

\def\nin{\noindent}
\def\cen{\centerline}
\def\pagebreak{\vskip 0pt plus 0.0001fil\break}
\def\linebreak{\break}

\def\hat{\widehat}

\def\eps{\varepsilon}
\def\epsilon{\varepsilon}
\def\eset{\emptyset}

\def\nin{\noindent}
\def\oline{\overline}

\def\pder#1,#2,#3 { {\partial #1 \over \partial #2}(#3)}
\def\pde#1,#2 { {\partial #1 \over \partial #2}}
\def\phi{\varphi}


\def\subeq{\subseteq}
\def\supeq{\supseteq}

\def\tilde{\widetilde}

\def\up{{\uparrow}}

\font\eightrm=cmr8


\font\smc=cmcsc10
\font\bfone=cmbx10 scaled\magstep1 
\font\bftwo=cmbx10 scaled\magstep2 

\def\qed{{\unskip\nobreak\hfil\penalty50\hskip .001pt \hbox{}\nobreak\hfil
          \vrule height 1.2ex width 1.1ex depth -.1ex
           \parfillskip=0pt\finalhyphendemerits=0\medbreak}\rm}

\def\qeddis{\eqno{\vrule height 1.2ex width 1.1ex depth -.1ex} $$
                   \medbreak\rm}

\def\Lemma #1. {\bigbreak\vskip-\parskip\noindent{\bf Lemma #1.}\quad\it}

\def\Sublemma #1. {\bigbreak\vskip-\parskip\noindent{\bf Sublemma #1.}\quad\it}

\def\Proposition #1. {\bigbreak\vskip-\parskip\noindent{\bf Proposition #1.}
\quad\it}

\def\Corollary #1. {\bigbreak\vskip-\parskip\nin{\bf Corollary #1.}
\quad\it}

\def\Theorem #1. {\bigbreak\vskip-\parskip\noindent{\bf Theorem #1.}
\quad\it}

\def\Definition #1. {\rm\bigbreak\vskip-\parskip\noindent{\bf Definition #1.}
\quad}

\def\Remark #1. {\rm\bigbreak\vskip-\parskip\noindent{\bf Remark #1.}\quad}

\def\Example #1. {\rm\bigbreak\vskip-\parskip\noindent{\bf Example #1.}\quad}

\def\Problems #1. {\bigbreak\vskip-\parskip\noindent{\bf Problems #1.}\quad}
\def\Problem #1. {\bigbreak\vskip-\parskip\noindent{\bf Problems #1.}\quad}

\def\Conjecture #1. {\bigbreak\vskip-\parskip\noindent{\bf Conjecture #1.}\quad}

\def\Proof#1.{\rm\par\ifdim\lastskip<\bigskipamount\removelastskip\fi\smallskip
            \noindent {\bf Proof.}\quad}

\def\Axiom #1. {\bigbreak\vskip-\parskip\noindent{\bf Axiom #1.}\quad\it}

\def\Satz #1. {\bigbreak\vskip-\parskip\noindent{\bf Satz #1.}\quad\it}

\def\Korollar #1. {\bbr\vskip-\parskip\nin{\bf Korollar #1.} \quad\it}

\def\Bemerkung #1. {\rm\bigbreak\vskip-\parskip\noindent{\bf Bemerkung #1.}
\quad}

\def\Beispiel #1. {\rm\bigbreak\vskip-\parskip\noindent{\bf Beispiel #1.}\quad}
\def\Aufgabe #1. {\rm\bigbreak\vskip-\parskip\noindent{\bf Aufgabe #1.}\quad}

\def\Beweis#1. {\rm\par\ifdim\lastskip<\bigskipamount\removelastskip\fi
           \smallskip\noindent {\bf Beweis.}\quad}

\nopagenumbers

\def\date{\ifcase\month\or January\or February \or March\or April\or May
\or June\or July\or August\or September\or October\or November
\or December\fi\space\number\day, \number\year}

\def\title{Title ??}
\def\author{Author ??}

\def\thanks#1{\footnote*{\eightrm#1}}

\def\rightheadline{\hfil{\eightrm\title}\hfil\tenbf\folio}
\def\leftheadline{\tenbf\folio\hfil{\eightrm\author}\hfil}
\headline={\vbox{\line{\ifodd\pageno\rightheadline\else\leftheadline\fi}}}

\def\firstheadline{}
\def\firstfootline{\cen{\rm\folio}}

\def\seite #1 {\pageno #1
               \headline={\ifnum\pageno=#1 \firstheadline
               \else\ifodd\pageno\rightheadline\else\leftheadline\fi\fi}
               \footline={\ifnum\pageno=#1 \firstfootline\else{}\fi}}

\newdimen\dimenone
 \def\checkleftspace#1#2#3#4{
 \dimenone=\pagetotal
 \advance\dimenone by -\pageshrink   
 \ifdim\dimenone>\pagegoal          
   \else\dimenone=\pagetotal
        \advance\dimenone by \pagestretch
        \ifdim\dimenone<\pagegoal
          \dimenone=\pagetotal
          \advance\dimenone by#1         
          \setbox0=\vbox{#2\parskip=0pt                
                     \hyphenpenalty=10000
                     \rightskip=0pt plus 5em
                     \noindent#3 \vskip#4}    
        \advance\dimenone by\ht0
        \advance\dimenone by 3\baselineskip   
        \ifdim\dimenone>\pagegoal\vfill\eject\fi
          \else\eject\fi\fi}


\def\subheadline #1{\nin\bigbreak\vskip-\lastskip
      \checkleftspace{0.7cm}{\bf}{#1}{\medskipamount}
          \indent\vskip0.7cm\centerline{\bf #1}\medskip}

\def\sectionheadline #1{\bigbreak\vskip-\lastskip
      \checkleftspace{1.1cm}{\bf}{#1}{\bigskipamount}
         \vbox{\vskip1.1cm}\cen{\bfone #1}\bsk}

\def\lsectionheadline #1 #2{\bigbreak\vskip-\lastskip
      \checkleftspace{1.1cm}{\bf}{#1}{\bigskipamount}
         \vbox{\vskip1.1cm}\cen{\bfone #1}\msk \cen{\bfone #2}\bsk}

\def\lchapterheadline #1 #2{\bigbreak\vskip-\lastskip\indent\vskip3cm
                       \cen{\bftwo #1} \msk \cen{\bftwo #2} \gsk}
\def\llsectionheadline #1 #2 #3{\bigbreak\vskip-\lastskip\indent\vskip1.8cm
\cen{\bfone #1} \msk \cen{\bfone #2} \msk \cen{\bfone #3} \nobreak\bsk\nobreak}


\newtoks\literat
\def\[#1 #2\par{\literat={#2\unskip.}%
\hbox{\vtop{\hsize=.15\hsize\nin [#1]\hfill}
\vtop{\hsize=.82\hsize\nin\the\literat}}\par
\vskip.3\baselineskip}

\mathchardef\emptyset="001F 
\def\address{Author: \tt$\backslash$def$\backslash$address$\{$??$\}$}

\def\firstpage{\nin
{\obeylines \parindent 0pt }
\vskip2cm
\centerline {\bfone \title}
\gsk
\centerline{\bf\author}

\vskip1.5cm \rm}

\def\lastpage{\par\vbox{\vskip1cm\nin
\line{
\vtop{\hsize=.5\hsize{\parindent=0pt\baselineskip=10pt\nin\address}}
\hfill} }}


\def\firstpage{\nin
{\obeylines \parindent 0pt }
\vskip2cm
\centerline {\bfone \title}
\ssk
\centerline {\bfone \titletwo}
\gsk
\centerline{\bf\author}
\vskip1.5cm \rm}

\def\bs{\backslash} 
\def\addots{\mathinner{\mkern1mu\raise1pt\vbox{\kern7pt\hbox{.}}\mkern2mu
\raise4pt\hbox{.}\mkern2mu\raise7pt\hbox{.}\mkern1mu}}

\pageno=1
\def\up#1{\leavevmode \raise.16ex\hbox{#1}}
 at 8truept
 at 8truept
 at 12truept
\chardef\ss="19
\def\3{\ss}
\def\bobheadline #1{\nin\bigbreak\vskip-\lastskip
      \checkleftspace{0.1cm}{\bf}{#1}{\medskipamount}
          \indent\vskip0.cm\centerline{\bf #1}\medskip}
\def\title{Holomorphic extensions of representations: (II) geometry}
\def\titletwo{and harmonic analysis}
\def\author{Bernhard Kr\"otz${}^*$ and Robert J. Stanton${}^\dagger$ }
\footnote{}{${}^*$ Supported in part by the NSF-grant DMS-0097314}
\footnote{}{${}^\dagger$ Supported in part by the NSF-grant DMS-0070742}

\def\date{November 25, 2001}
\def\Box #1 { \msk\par\nin 
\centerline{
\vbox{\offinterlineskip
\hrule
\hbox{\vrule\strut\hskip1ex\hfil{\smc#1}\hfill\hskip1ex}
\hrule}\vrule}\msk }

\def\bs{\backslash} 

\def\address
{The Ohio State University 

Department of Mathematics 

231 West 18th Avenue 

Columbus, OH 43210--1174 

USA

{\tt kroetz@kurims.kyoto-u.ac.jp}

{\tt stanton@math.ohio-state.edu}

}

\firstpage 

\bobheadline{Introduction}

\par Akhiezer and Gindikin [AG90] proposed the existence of a natural complexification of a Riemannian symmetric space of noncompact type, $G/K$.  They sought a complex $G$-space with proper $G$ action and which contained $G/K$ as a totally real embedded submanifold. The obvious candidate for a complexification, $G_\C/K_\C$, fails to have the $G$ action proper, nevertheless in it they identified a candidate, $$\Xi \subeq G_\C/K_\C.$$ 
The open subset $\Xi$ inherits a complex structure from $G_\C/K_\C$,  and it contains $G/K$ as a totally real submanifold via the natural inclusion $G/K\hookrightarrow G_\C/K_\C$.

\par Almost simultaneously a related development was taking place in complex differential geometry. In the 60's much work had been done determining which geometric structures on a manifold $M$ prolong to similar geometric structures on $TM$. In the 90's a significant refinement of Grauert's embedding theorem was developed. A real analytic, compact manifold with a complete Riemannian metric, $(M,g)$, was shown to prolong to further structure. Namely, one obtains on an open neighborhood of the zero section in $TM$, a so-called Grauert domain, a unique complex structure having the property that the canonical lift of any geodesic to a map from $T\Bbb R\cong \Bbb C \to TM$ is holomorphic. Using the metric one can identify the cotangent bundle with the tangent bundle. If one combines this complex structure with the transported symplectic structure, one obtains a K\"ahler metric on the Grauert domain.

\par The Akhiezer-Gindikin domain can easily be shown to be a Grauert domain for the Riemannian symmetric space $G/K$. One of the first results in this paper, a consequence of Theorem 2.4, is that $\Xi$ is a $G$-invariant domain of holomorphy. The proof involves a major improvement of one of the main results from [KS04]. Briefly, there we identified a natural complexification for semisimple groups $G\subset  G_\C $ as an open $G-K_\C$ double coset of the form $G\exp(i\Omega) K_\C$. The importance of this domain is that we showed that any $K$-finite matrix coefficient of an irreducible unitary representation of $G$ has holomorphic extension to it, and blows-up at specific points of the boundary of $\Omega$.  Fortuitously, this construction and that of Akhiezer-Gindikin are compatible in that $G\exp(i\Omega) K_\C/K_\C$ and $ \Xi$ are biholomorphic. Then to obtain the result, we show in Theorem 2.4 that for spherical unitary representations the holomorphic extension blows-up at!
  all points of the boundary of $\Xi$.  
\par Hence any symmetric space of noncompact type $X =  G/K$ prolongs to a canonical  complex manifold $\Xi$ contained in the tangent bundle $TX$, and $\Xi$ is a natural domain of holomorphy for harmonic analysis on $G/K$. One of the goals in the paper is to formulate some fundamental harmonic analysis on $G/K$ in terms of complex analysis on $\Xi\subset TX$. The other main goal of the paper is to give a detailed description of the complex differential geometry arising from several natural K\"ahler structures defined on $\Xi$ or natural subdomains of it. Besides the K\"ahler structure associated to the canonical symplectic structure on the cotangent bundle, there are two more surprising K\"ahler structures that we identified and whose associated metric holomorphically extends that of $G/K$. 
 \par  For applications to harmonic analysis one needs the detailed structure of the Riemannian metric on $\Xi$, the Riemannian measure on $\Xi$, etc. This we do in the first part of the paper for the canonical K\"ahler structure on $\Xi$. In \S 2 we show that $\Xi$ has an ample supply of plurisubharmonic exhaustion functions, and from this obtain an easy proof that $\Xi$ is Stein. While one always has existence of the canonical K\"ahler structure on the Grauert domain, one seldom has an explicit formula for it. In \S 4 we describe in Lie theoretic terms the canonical K\"ahler structure on $\Xi$. Similarly for the associated Riemannian structure, we present an explicit expression for the Riemannian measure in terms of roots. For this Riemannian structure we identify the element in the universal enveloping algebra that induces the Laplace-Beltrami operator on $\Xi$. The results are consistently presented from the point of view of possible applications to harmonic analysis, i.!
 e. using Lie structure theory. 
\par The reformulation of harmonic analysis in complex analytic terms originated in [KS04] and was presented there in completeness for classical groups. In \S 1 this is extended to exceptional groups. The holomorphic extension was critical for the recent successes on the Gelfand-Gindikin program in [GK\'O03] and [GK\'O04]. Here, in a direction very different from that topic, we present new connections of representation theory with complex differential geometry.  For example, we show in \S 3 that $\hat G_s$, the $K$-spherical irreducible unitary dual of $G$, embeds in a specific moduli space of K\"ahler structures on $\Xi$. While tempted to conjecture that the spherical unitary dual is parametrized by such K\"ahler metrics, the natural approach to this problem would require information about solutions with prescribed singularities to a Monge-Amp\`ere equation on a noncompact manifold. Further applications to harmonic analysis are presented in \S 5 and \S6. The use of Hilbert !
 spaces of holomorphic functions is a familiar technique in the study of highest weight representations. In \S 5 we show that $K$-spherical unitary representations are similarly realizable by Hilbert spaces of holomorphic functions - a result new even for $Sl(2,\Bbb R)$.  
The heat kernel associated to the Riemannian Laplace operator is fundamental to harmonic analysis. In  \S 6 we give a holomorphic extension of the heat kernel on $G/K$ to $\Xi$.  One immediate application is to define a $G$-equivariant holomorphic heat kernel transform (i.e. Bargmann-Segal transform) for $G/K$. 
 \par Whereas there is only one complex structure on $\Xi$ that we should consider, there are several K\"ahler metrics whose associated Riemannian metric will extend the metric on $G/K$. These other K\"ahler metrics need not be defined on all of $\Xi$.  In \S 7 we define a subdomain, denoted $\Xi_0$, that we show is 
bi-holomorphic to a Hermitian symmetric space but with the inherited metric not isometric to it.  We also identify which domains have $\Xi = \Xi_0$. For these $\Xi$ one has then an unexpected hidden symmetry group  not predictible from their structure theoretic definition, and one can impose on $\Xi$ a Hermitian symmetric metric which also extends the metric on $G/K$. That eigenfunctions of the 
invariant differential operators on $G/K$ and, more generally, quotients of $G/K$ by lattices in $G$ have holomorphic extension to a Hermitian symmetric space seems very intriguing.  Finally, in \S 8 we describe and characterize a subdomain $\Xi^{1,{1\over2}}$ of $\Xi_0$ and produce a canonical K\"ahler metric on it using Jordan algebra formulas. This metric too gives an extension of the Riemannian metric on $G/K$ but to a different subdomain of $\Xi$. What role this domain and metric  has in harmonic analysis or automorphic functions is not at all clear to us.
\vfill 
\eject
\par Since our first lecture on our results in April 2000 at Oberwolfach, and since the original submission of this manuscript in February 2002, various individuals obtained, frequently with different methods, some results presented here. We cite the work of which we are aware at the appropriate point in this paper.

\par We want to express our sincere gratitude to L. Barchini, S. Gindikin, 
K.-H. Neeb  and G. \'Olafsson for several useful conversations while we completed this work.

\bobheadline{\S 1 Holomorphic extensions}

In [KS04] we gave a proof of the maximal holomorphic
extension of the orbit map of a $K$-finite vector of
an irreducible Hilbert representation of a connected classical semisimple Lie group.  
Conjecture A in [KS04] was the technical result needed to extend the proof to 
all connected semisimple Lie groups. Now from results in 
[Hu02] or [M03], see also a related result in [B03], this is known to be proved. 
For the readers convenience and because of its importance to the paper we sketch how   the holomorphic extension follows from it.

Let $G$ be a connected reductive Lie group contained in a complexification
$G_\C$.  We denote by $\g$ (resp.  $\g_\C$) the Lie algebra of $G$ (resp.
$G_\C$).  Let $K$ be a maximal compact subgroup of $G$ and let $\k$ be its
Lie algebra.  Attached to $\k$ and the Killing form is the Cartan
decomposition $\g=\k\oplus\p$.  Take $\a\subeq \p$ a maximal abelian subspace
and let $\Sigma =\Sigma (\g, \a)\subeq \a^*$ be the corresponding system of
roots.  The simultaneous eigenspaces of $\ad (H), H\in \a$, give the root
space decomposition

$$\g=\a\oplus \m \oplus\bigoplus_{\alpha\in \Sigma} \g^\alpha,$$ here
$\m=\z_\k(\a)$ and $\g^\alpha =\{ X\in \g\:  (\forall H\in \a) [H,
X]=\alpha(H) X\}$.  For each choice of positive roots,
$\Sigma^{+}\subeq\Sigma$, one obtains a nilpotent Lie algebra
$\n=\bigoplus_{\alpha\in\Sigma^+}\g^\alpha$ and an
Iwasawa decomposition on the Lie algebra level 
$$\g=\k\oplus \a \oplus \n.$$ Denote by $A$ and $N$
the analytic subgroups of $G$ corresponding to $\a$ and $\n$.  These choices
give an Iwasawa decomposition for $G$; namely, the multiplication map
$$N\times A\times K\to G, \ \ (n, a, k)\mapsto nak$$ is an analytic
diffeomorphism.  Thus, every element $g\in G$ can be written uniquely as
$g=n(g)a(g)\kappa(g)$, with each of the maps $n(g)\in N$, $a(g)\in A$,
$\kappa(g)\in K$ depending analytically on $g\in G$.  The last piece of 
structure
theory we shall recall is the Weyl group of $\Sigma(\g,\a)$, ${\cal W}
=N_K(\a)/ Z_K(\a)$.

With Akhiezer-Gindikin we define a domain in $G_\Bbb C$ using the restricted
roots.  Set $$\Omega=\{ X\in\a\:(\forall \alpha\in \Sigma)\
|\alpha(X)|<{\pi\over2}\}.$$ Then $\Omega$ is convex and ${\cal
W}$-invariant.  Set

$$T_\Omega=A \exp(i\Omega).$$ 
$T_\Omega$ is open in $A_\C$ and diffeomorphic to 
$ A\times \Omega$, which, in turn, is diffeomorphic to $ \a \oplus
i\Omega$. In other words $T_\Omega$ is diffeomorphic to a tube domain in
a complex vector space.

We showed in [KS04] that the domain $G\exp(i\Omega)K_\C$ is an open
$G$-$K_\C$ double coset in $G_\C$. For classical groups we also showed that
$G\exp(i\Omega) K_\C \subeq N_\C A_\C K_\C $, and that the Iwasawa projection,
$a\:  G\to A$, extends holomorphically to a map $a\:  G\exp(i\Omega) K_\C \to
A_\C $ such that $x\in N_\C a(x) K_\C$.  Now, as a consequence of the
complex convexity theorem of Gindikin-Kr\"otz (cf.\ [GK02b]), we have more 
precisely that

$$G\exp(i\Omega) K_\C \subeq N_\C T_\Omega K_\C, \ \leqno(1.1)$$  giving a
holomorphic extension of the Iwasawa projection $a\: G\exp(i\Omega) K_\C \to
A_\C $ with values in $T_\Omega$. Exactly as in
[KS04, Th.\ 3.1] one obtains 

\Theorem 1.1.  Let $(\pi, {\cal H})$ be an admissible Hilbert representation
of $G$.  Then for any $K$-finite vector $v\in {\cal H}$, the orbit map $$G\to
{\cal H}, \ g\mapsto \pi(g)v$$ extends to a $G$-equivariant holomorphic map
on $G \exp(i\Omega)K_\C$.\qed

>From Theorem 1.1 follows immediately

\Corollary 1.2.  Let $(\pi,{\cal H})$ be an admissible Hilbert representation
of $G$.  Then for all $K$-finite vectors $v\in {\cal H}$ and for all
hyperfunction vectors $\eta\in {\cal H}^{-\omega}$, the generalized matrix
coefficient $$\pi_{v,\eta}\:  G\to \C, \ \ g\mapsto
\la\pi(g)v,\eta\ra \:=\oline {\eta(\pi(g)v)}$$ extends to a holomorphic
function on $G\exp(i\Omega)K_\C$.  \qed

If one takes right $K_\C$ cosets of $G\exp(i\Omega)K_\C$, one obtains
the domain $$\Xi\:= G\exp(i\Omega)K_\C/ K_\C, $$ introduced in [AG90]
and dubbed by Gindikin the {\it complex crown} of the Riemannian symmetric
space $G/K$. 
\par  $\Xi$, as an open subset of $G_\C /K_\C$, inherits the 
structure of a complex manifold. In this paper we will obtain several facts about the complex structure as 
well as geometric properties of this 
domain and domains related to it. Also,  
the relationship of harmonic analysis on $G/K$ and complex analysis on $\Xi$ will 
be investigated. The following is crucial to these developments. 
\par Let $\D(G/K)$ denote the algebra of $G$-invariant
differential operators on $G/K$.  Specializing Corollary 1.2 to $K$-fixed
vectors, and using the solution by [K-T78] to the Helgason conjecture, one
obtains 

\Proposition 1.3.  Every eigenfunction $\phi$ on
$G/K$ of the algebra of $G$-invariant differential operators $\D(G/K)$
extends holomorphically to $\Xi$.\qed

\bobheadline{\S 2 Plurisubharmonic functions on $\Xi$}
\par From Proposition 1.3, for $(\pi_\lambda,{\cal H}_\lambda)$ a unitary spherical principal series representation with parameter $\lambda\in\a_\C^\ast$ and with 
non-zero $K$-fixed vector $v_0$, it follows that the 
spherical function $\phi_\lambda$ extends to a  holomorphic function on $\Xi$. In [KS04] we showed that this holomorphic extension blows-up in specific directions at the boundary. In this section we shall show how the holomorphic extension of the spherical function $\phi_\lambda$ provides $\Xi$ blows-up at all points of the boundary and provides $\Xi$ with a supply of strictly pluriharmonic exhaustion functions.

\par For $Z\in\g_\C$ we define a vector field $\tilde Z$ on $\Xi$ by
$$(\tilde Z f)(x)={d\over dt}\Big|_{t=0} f(\exp(tZ)x),$$ here $f$ is a function 
on
$\Xi$ differentiable at $x$. Let $J$ denote the linear operator on $\g_\C$ 
corresponding to multiplication by $i$ on $\g_\C$. To $Z$ one can associate the 
familiar Cauchy-Riemann operators $\partial_Z\:={1\over 2} (\tilde Z -i 
\tilde{JZ})$ and 
$\oline\partial_Z\:={1\over 2} (\tilde Z +i \tilde{JZ})$.
\par Our first observation relating the unitary dual to complex analysis on 
$\Xi$ is 
given in the following result.

\Proposition 2.1.  Let $(\pi, {\cal H})$ be a non-trivial $K$-spherical unitary 
irreducible representation of $G$. Let $v_0$ be a normalized  $K$-spherical vector. 
Then the functions 
$$s_\pi\:\Xi\to \R, \ \ xK_\C\mapsto \|\pi(x)v_0\|^2\ , $$ 
and 
$$\log (s_\pi)\:\Xi\to \R, \ \ xK_\C\mapsto 2\log \|\pi(x)v_0\|$$
are  $G$-invariant and 
strictly plurisubharmonic.

\Proof. Before we start the proof, let us recall some standard 
facts on plurisubharmonic functions. Let  $M$ be  a complex
manifold, ${\cal H}$ a Hilbert space and $f\: M\to {\cal H}$
a holomorphic map. Then both 
$$z\mapsto \|f(z)\|^2, \qquad z\mapsto \log \|f(z)\|$$
are plurisubharmonic functions on $M$. The
first one is strictly plurisubharmonic if and only if  $f$ is an immersion.
The second map is strictly plurisubharmonic if and only if for all $m\in M$
the vector $df(m)(Z)$ is transversal to $f(m)$ for all $Z\neq 0$ in a fixed 
totally real subspace of $T_m M$. 

\ssk From Theorem 1.1 the orbit map from $G$ to ${\cal H}$, $g\mapsto 
\pi(g)v_0$, admits a holomorphic extension to a $G$-equivariant map from 
$G\exp(i\Omega)K_\C$ to ${\cal H}$, $\ x\mapsto \pi(x)v_0$.  Clearly this map is 
right 
$K_\C$-invariant, so factors to a holomorphic $G$-equivariant map 
$$f\: \Xi\to {\cal H}, \ \ xK_\C\mapsto \pi(x)v_0.$$  
The $G$-invariance of $s_\pi$ and $\log (s_\pi)$  follows from the 
unitarity of $\pi$.
\par In view of our comments at the beginning, the proof of the proposition 
will be complete once we have established the transversality of $f$. 
For that, fix $m=xK_\C \in \Xi$ and identify the tangent space 
$T_m\Xi$ with $\Ad(x)\p_\C$. Then $\Ad(x)\p$ becomes a totally real 
subspace of $T_m\Xi$. We compute $df(m)(Z)=d\pi(Z)\pi(x)v_0$. For 
$Z=\Ad(x)Y$ with $Y\in \p$, we thus get $df(m)(Z)=\pi(x)d\pi(Y)v_0$. 
Assume that $f$ is not transversal. Then for some $Y\in \p$, $Y\neq 0$, we would have 
$\pi(x)d\pi(Y)v_0\in \C\pi(x)v_0$ or, equivalently,  
$d\pi(Y)v_0\in \C v_0$. 
Now $\l\:=\{X\in \g\:  d\pi(X)v_0\in \C v_0\}$ is a subalgebra of $\g$ and it contains $\k$ and  $Y\in \p$. But then $\l=\g$, and so $\pi$ is the trivial representation - a 
contradiction. So $f$ is transversal. 
\qed

We continue with some properties of the function $\log (s_\pi)$. 

\Proposition  2.2. Let $(\pi, {\cal H})$ be a non-trivial $K$-spherical unitary 
irreducible representation 
of $G$. Then the function 
\item {(i)}$\log (s_\pi)\: \Xi\to \R$ is $G$-invariant and strictly
plurisubharmonic; 
\item{(ii)}$\log(s_\pi)\geq 0$;
\item{(iii)} $G/K=\{z\in \Xi\: \log(s_\pi)(z)=0\}$.

\Proof. (i) was already established in  Proposition 2.1. 
\par Moving on to (ii) and (iii). As $\log(s_\pi)$ is $G$-invariant it is 
specified by its values on $T_\Omega$. But by $G$-invariance and ${\cal W}$- 
invariance of $\log(s_\pi)$ we may consider $\log(s_\pi)$ as an $\{A, {\cal 
W}\}$-invariant function on the abelian tube domain 
$T_\Omega=A\exp(i\Omega)$. Then $\psi_\pi=\log(s_\pi)\circ \exp\res_\Omega$
defines a strictly convex ${\cal W}$-invariant function on $\Omega$. 
Clearly, (ii) and (iii) will be proved if we can show that $\psi_\pi$ has a
unique absolute minimum at $0$. Let us first show that $0$ is a minimum of 
$\psi_\pi$. Since 
$\psi_\pi=\psi_\pi\circ w$ for all $w\in {\cal W}$, and $w(0)=0$, we have 
 $$d\psi_\pi(0)=d\psi_\pi(0)\circ w$$
for all $w\in {\cal W}$. Hence $d\psi_\pi(0)=0$ and $0$ is a minimum, as 
$\psi_\pi$ is strictly convex. To conclude the proof it suffices to show that 
$0$ is the only minimum of $\psi_\pi$. If not, then there
exists an $X_0\in \Omega$, $X_0\neq 0$ such that $\psi_\pi$ has a local minimum 
at 
$X_0$, say $\psi_\pi(X_0)=m_0$. As $X_0\neq 0$ we find a $w\in {\cal W}$ such 
that 
$X_1=w(X_0)\neq X_0$. Since $\psi_\pi(X_1)=\psi_\pi(X_0)=m_0$, 
restricting $\psi_\pi$ to the line segment connecting $X_0$ and $X_1$ we 
obtain a strictly convex function with local minima at the endpoints - 
impossible. Hence there is no minimum other than $0$. \qed 

We denote the topological boundary of $\Xi$ in $G_\C/K_\C$ by $\partial 
\Xi$. We wish to investigate  the boundary behaviour of the plurisubharmonic 
functions $s_\pi$ and $\log(s_\pi)$. We start with a geometrical fact.

\Lemma 2.3. \item{(i)}
$G\exp(i\partial \Omega)K_\C/ K_\C\subeq \partial\Xi$.  \item{(ii)} Let
$\{x_n=g_n\exp(iX_n)K_\C \: g_n\in G, X_n\in \Omega \}$ be a
sequence in $\Xi$ with $x_n\to x\in \partial\Xi$.  Then $X_n\to X\in \partial
\Omega$.

\Proof.  (i) This is proved in [AG90, p.\ 8-9].  \par\nin (ii) Let $\theta$
also denote the complex linear extension to $G_\C$ of the Cartan involution on 
$G$.  Recall from (1.1) that $G\exp(i\Omega)\subeq N_\C A_\C K_\C$. Then 
$g_n\exp(iX_n)\in N_\C A_\C K_\C$ and
$$y_n =g_n\exp(i2X_n)\theta(g_n)^{-1}\in N_\C A_\C \oline N_\C,$$ where,as 
usual,
$\oline N_\C =\theta(N_\C)$.  Write $y\in G_\C$ for the limit of
$(y_n)_{n\in \N}$.  \par Suppose that there is a subsequence of $(X_n)_{n\in 
\N}$ which converges to an element in $\Omega$.
Replacing the original sequence by the corresponding subsequence, we may assume 
that$X_n\to X\in \Omega$.  
Write $g_n =b_n k_n$ for $b_n\in AN$ and $k_n\in K$. As $K$ is compact, we may 
assume that $k_n\to k\in K$.  
Set $z_n =k_n\exp(i2X_n)k_n^{-1}$. Then $\{z_n\}$ converges, say to 
$z=k\exp(i2X)k^{-1}$, and it follows from (1.1) 
that $z\in N_\C A_\C \oline N_\C$.  If $(b_n)_{n\in \N}$ were unbounded, then 
from [KS04, Lemma 1.6] 
we would have $(y_n)_{n\in \N}$ an unbounded sequence in $G_\C$ - a 
contradiction, as $y_n\to y$.  Hence$(b_n)_{n\in \N}$ is bounded.  
Thus we may assume that $(g_n)_{n\in \N}$ is
convergent to some element $g\in G$.  But then $x$ must be in $\Xi$, 
contradicting $x\in
\partial \Xi$.  \qed

We come now to the main result in this section.

\Theorem 2.4.  Let $(\pi, {\cal H})$ be a non-trivial $K$-spherical unitary 
irreducible
representation of $G$.  Then there exists a constant $C_\pi>0$ such that 
$$(\forall X\in \partial \Omega)\qquad s_\pi(\exp(i(1-\eps)X)K_\C)\geq 
C_\pi|\log \eps|\leqno(2.1)$$ for all $0<\eps\leq 1$.  In particular, one has 
$$\lim_{x\to \partial \Xi\atop x\in\Xi} s_\pi(x)=\infty.$$

\Proof.  The last assertion of the theorem follows from (2.1) and
Lemma 2.3(ii).  So we are left with proving (2.1).  \par Let $X\in
\partial\Omega$.  It follows from the definition of $\Omega$ that there
exists an $\alpha\in\Sigma$ such that $\alpha(X)={\pi\over 2}$.  Write
$H_\alpha\in\a$ for the co-root of $\alpha$.  Take
$X_\alpha\in\g^\alpha$ such that $[X_\alpha, \theta X_\alpha]= -H_\alpha$. Then  
$$\g_0 =\span_\R\{H_\alpha, X_\alpha, \theta X_\alpha\}$$ is a
$\theta$-stable subalgebra of $\g$ isomorphic to $\sL(2,\R)$.
Let $\a_\alpha =\alpha^\bot$ be the kernel of $\alpha$.  Notice that $$\g_1 \: 
=\g_0 + \a_\alpha$$ is in fact a Lie algebra direct sum.  Denote by $G_1$, $G_0$ 
and $A_\alpha$ the analytic
subgroups of $G$ with Lie algebras $\g_1$, $\g_0$ and $\a_\alpha$. One has
$$G_1=G_0\times A_\alpha.$$ Now if ${\pi\over 4}H_\alpha\in
\partial\Omega$, then a mild generalization of Proposition 4.5 in [KS04] 
(simply replace $G_0$ by $G_1$ in the proof) gives

$$s_\pi(\exp(i(1-\eps)X))\geq C_\alpha |\log\eps|,$$ for all $0<\eps\leq
1$ and a constant $C_\alpha>0$ depending only on $\alpha$.  Thus it is enough to 
show:  \msk \item{(2.2)} Among all $\alpha\in
\Sigma$ with $\alpha(X)={\pi\over 2}$ there is one with ${\pi\over
4}H_\alpha\in \partial \Omega$.  \msk\nin We will verify (2.2) by considering
the different types of root systems.  Of course we may assume that $\g$ is
simple.  (2.2) would be a consequence of the following statement (2.3),  
however, (2.3) is not valid for all root systems. \msk
\item{(2.3)} For all $\alpha\in \Sigma$, ${\pi\over
4}H_\alpha\in\partial\Omega$.
\msk\nin We proceed to consider the various root systems.
\ssk\nin Case 1:  $\Sigma$ is of type $A_n$, $C_n$ or $D_n$.  \par\nin In all
these cases (2.3) holds (cf.\ [KS04, Sect.\ 4]), and thus so does (2.2).  
\ssk\nin Case 2:  $\Sigma$ is of type $BC_n$.
\par\nin This can be reduced to the $C_n$ case since any root
$\alpha\in\Sigma(BC_n)$ with $\alpha(X)={\pi\over 2}$ must lie in the natural 
$\Sigma(C_n)\subeq
\Sigma(BC_n)$.  \ssk\nin Case 3:  $\Sigma$ is of type $B_n$.  \par\nin Recall 
that
the roots of $B_n$, in the standard notation, are $$\Sigma(B_n)=\{\pm \eps_i
\pm\eps_j\:  1\leq i\neq j\leq n\}\cup\{ \pm \eps_i\:  1\leq i\leq n\},$$ and
for a basis of $\Sigma (B_n)$ we may take
$$\Pi=\{ \alpha_1, \ldots, \alpha_n\}=\{ \eps_1-\eps_2, \eps_2-\eps_3,
\ldots, \eps_{n-1}-\eps_n, \eps_n\}.$$ Define $e_j\in \a$ by
$\eps_j(e_i)=\delta_{ij}$. Take $X\in \partial\Omega$.  By the Weyl
group invariance of (2.1) we may assume that $\alpha_i(X)={\pi\over 2}$ for
some $1\leq i\leq n$.  Now ${\pi\over 4}H_{\alpha_i}\in \partial\Omega$ for
$1\leq i\leq n-1$, while ${\pi\over 4}H_{\alpha_n}\not\in\partial\Omega$.
Thus we have to consider only the case where $\alpha_n(X)={\pi\over 2}$.  But
this determines $X$ uniquely to be $X={\pi\over 2} e_n$ (see the proof of
[KS04, Lemma 2.8]).  Hence $-\alpha_{n-1}(X)={\pi\over 2}$, and (2.2) is
established in this case.

\ssk\nin Case 4:  $\Sigma$ is of type $E_6$, $E_7$ or $E_8$.  \par\nin We
will show that (2.3) holds for all these cases.  We have to consider only the
case of $\Sigma(E_8)$ since $\Sigma (E_6)\subeq \Sigma (E_7)\subeq \Sigma 
(E_8)$.  Now the roots of $E_8$ are 
$$\Sigma(E_8)=\{ \pm\eps_i\pm\eps_j\:  1\leq i\neq j\leq n\}\cup \{{1\over 2}
\sum_{j=1}^8 (-1)^{n(j)} \eps_j\:  \sum_{j=1}^n n(j)\ \hbox{even}\},$$
and for a basis of $\Sigma(E_8)$ we may take

$$\eqalign{\Pi&=\{ \alpha_1,\ldots, \alpha_8\}= \{{1\over
2}(\eps_8-\eps_7-\eps_6-\eps_5-\eps_4-\eps_3-\eps_2+\eps_1),\cr
&\eps_2+\eps_1, \eps_2-\eps_1, \eps_3-\eps_2, \eps_4-\eps_3, \eps_5-\eps_4,
\eps_6-\eps_5, \eps_7-\eps_6\}.\cr}$$ In order to establish (2.3) it is
enough to show that ${\pi\over 4}H_{\alpha_i}$ is in $\partial\Omega$
for all $1\leq i\leq 8$.  Now $\alpha_2, \ldots, \alpha_7$ span
a root subsystem $\Sigma_0$ of $\Sigma(E_8)$ of type $D_7$, and (2.3) holds for 
such a system.  From this it is easy
to see that ${\pi\over 4}H_\alpha\in \partial\Omega$ for all $\alpha\in
\Sigma_0$.  It remains to consider $$H_{\alpha_1}={1\over
2}(e_8-e_7-e_6-e_5-e_4-e_3-e_2+e_1).$$ Here one readily checks that ${\pi\over
4}H_{\alpha_1}\in \partial\Omega$.

\ssk\nin Case 5:  $\Sigma$ is of type $G_2$.  \par\nin With
$$\a =\{\sum_{j=1}^3 x_j e_j\:  x_1+x_2+x_3=0\},$$ we have

$$\Sigma(G_2)=\{\pm(\eps_1-\eps_2),\pm(\eps_1-\eps_3) \pm(\eps_2-\eps_3)\}
\cup\{ \pm(2\eps_1 -\eps_2 -\eps_3),\pm(2\eps_2 -\eps_1 -\eps_3) \pm(2\eps_3
-\eps_1 -\eps_2)\}.$$ For a basis of $\Sigma(G_2)$ we take
$$\Pi=\{\alpha_1, \alpha_2\}=\{ \eps_1-\eps_2, -2\eps_1 +\eps_2+\eps_3\}.$$ Let 
$X\in\partial\Omega$.  Arguing as in Case 3 we may assume that
$\alpha_1(X)={\pi\over 2}$ or $\alpha_2(X)={\pi\over 2}$.  First we show
that $\alpha_1(X)={\pi\over 2}$ is impossible.  Write
$X=\sum_{j=1}^3 x_j e_j$ with $\sum_{j=1}^3 x_j=0$. Now if 
$\alpha_1(X)={\pi\over 2}$, we have $X={\pi\over 4}(e_1-e_2)$ - but then $X\in 
\partial\Omega$. This is not possible since
$$(2\eps_1 -\eps_2 -\eps_3)(X)>{\pi\over 2}.$$ Suppose  
$\alpha_2(X)={\pi\over 2}$. From the explicit description of $\Sigma (G_2)$ we 
see that $${\pi\over 4}H_{\alpha_2}={\pi\over 4} {1\over 2}(-2e_1 +e_2
+e_3) \in \partial\Omega, $$ thus establishing (2.2) for $\Sigma(G_2)$.

\ssk\nin Case 6:  $\Sigma$ is of type $F_4$.  \par\nin Here we have
$$\Sigma(F_4)=\{ \pm \eps_i\pm \eps_j\:  1\leq i\neq j\leq 4\}\cup
\{\pm\eps_i\:1\leq i\leq 4\} \cup\{ {1\over 2}(\pm \eps_1\pm
\eps_2\pm\eps_3\pm\eps_4)\},$$ and for a basis we take $$\Pi=\{\alpha_1, 
\alpha_2, \alpha_3,
\alpha_4\}= \{{1\over 2}(\eps_1-\eps_2-\eps_3-\eps_4), \eps_4, \eps_3-\eps_4,
\eps_2-\eps_3\}.$$ As explained in Case 3 we may assume that
$\alpha_i(X)={\pi\over2}$ for some $1\leq i\leq 4$.  Now it is easy to check
that ${\pi\over 4}H_{\alpha_j}\in \partial\Omega$ for $j=1,3,4$ but
${\pi\over 4}H_{\alpha_2}\not\in \partial\Omega$.  Thus the crucial case is
when $\alpha_2(X)={\pi\over 2}$.  But this determines $X$ uniquely, namely
$X={\pi\over 2}e_4$.  Then $-\alpha_3(X)={\pi\over 2}$ and (2.2) is proved
for the case of $\Sigma(F_4)$.\qed

\par Recall that $\Xi$ is open in $G_\C / K_\C$, so it inherits a complex 
structure. For $\pi_{\lambda}$ a $K$-spherical representation let
$\phi_{\lambda}(gK)=\la \pi_{\lambda}(g)v_0, v_0\ra\qquad  (g\in G)$
be the associated spherical function. In [KS04, \S 4] we showed that this function extends holomorphically 
to $\Xi$ and its restriction to $A$ extends holomorphically 
to $T(2\Omega)$. Denote these holomorphic extensions also 
by $ \phi_{\lambda}$. The relation with $s_{\pi_{\lambda}}$ is 
$(\forall a\in \exp(i\Omega))\qquad \phi_\pi(a^2)=s_\pi(a).$ Then from Theorem 2.4 it follows that $\Xi$ is a $G$-invariant domain of holomorphy. See also [F03] and [B03].
\par Later we explain why this complex structure is the adapted complex 
structure arising from the Riemannian structure on $G/K$, but for now we are 
concerned with properties related to the homogeneous nature of $\Xi$. One knows 
that $G$ acts properly on $\Xi$ ([AG90, Prop.\ 3]). So if
$\Gamma$ is a torsion-free discrete subgroup of $G$, the quotient
$\Gamma\bs \Xi$ is a complex manifold.

As an immediate consequence of Theorem 2.4 we have

\Proposition 2.5.  Let $\Gamma<G$ be a torsion free co-compact subgroup.  Then
the complex manifold $\Gamma\bs \Xi$ is Stein.

\Proof.  Let $s_\pi\:  \Xi\to \R$ be the function from Proposition 2.1.  Since
$s_\pi$ is left $G$-invariant, this function factors to a strictly
plurisubharmonic function on $\Gamma\bs \Xi$ denoted $s_\pi^{\Gamma}$.  In order
to show that $\Gamma\bs \Xi$ is Stein, it suffices to show that $s_\pi^{\Gamma}$ 
is
proper.  For that let $(x_n)_{n\in \N}$ be a sequence which tends to infinity
in $\Gamma\bs \Xi$.  Then $x_n = \Gamma g_n \exp (iY_n) K_\C$ and, as $\Gamma\bs 
G$ is compact, $x_n\to\partial\Xi$. Then Theorem 2.4 implies   
$$\lim_{n\to\infty}s_\pi^{\Gamma}(x_n)=\lim_{n\to\infty}s_\pi(\exp(iY_n))= \infty 
.$$ \qed

\Proposition 2.6. $\Xi$, with its inherited complex structure from $G_\C /K_\C$, 
is Stein.

\Proof.  We have only to construct a strictly plurisubharmonic exhaustion
function on $\Xi$.  Since $G_\C/ K_\C$ is Stein, there exists a
plurisubharmonic exhaustion function on $G_\C/K_\C$, $\psi\:  G_\C/ K_\C \to 
\R^+$.  Now define a function $\phi \:  \Xi\to \R^+$ by $$\phi =\psi\res_\Xi+
s_\pi,$$ where $s_\pi$ is a function as in Proposition 2.1.  As $\phi$ is a sum 
of two strictly plurisubharmonic functions, $\phi$ is strictly plurisubharmonic. 
 It remains to show that $\phi$ is proper. Suppose that $x_n\to \infty$.  That 
$x_n\to \infty$ means either $x_n\to x\in \partial \Xi$ or $x_n\to\infty$ in
$G_\C/ K_\C$. As both $\psi$ and $ s_\pi$ are positive, in the first case it 
follows from (2.1) that $\lim_{n\to
\infty}\phi(x_n)=\infty$, while in the latter case it follows from the fact that 
$\psi$ is an exhaustion function.\qed

\Remark 2.7. The result presented in Proposition 2.6 has been approached from 
several directions. We refer the reader to [AG90], [B03], [BHH03], [GK02b], [GM01] 
and [Hu02].

\bobheadline{\S 3  A moduli space based on $\Xi$ }
\par We shall denote by $\hat G_s$ the $K$-spherical unitary 
dual of $G$. In this section we show that $\hat G_s$ provides 
the complex manifold $\Xi$ with a space of K\"ahler potentials.
\par For $\pi\in \hat G_s-\{ 1\}$ we write $g^\pi$ for the 
Riemannian metric associated to the $G$-invariant strictly 
plurisubharmonic function $\log(s_\pi)$. Explicitly it is 
given by 

$$g_z^\pi(\tilde Z_z, \tilde Y_z)=\Re \left(\partial_Z \oline \partial_Y \log(s_\pi)\right)(z) \qquad 
(z\in \Xi) (Z,Y\in\g_\C).$$
Below we will show that $g^\pi$ is complete. For that 
it is useful to recall some facts on the infinite dimensional projective space. 

\par Let ${\cal H}$ denote a complex Hilbert space and let 
$\P({\cal H})$ be its associated projective space. For 
$v\in {\cal H}\bs \{0\}$ we write $[v]=\C v$ for the corresponding 
line through $v$. Then there is the projection mapping 

$$p\: {\cal H}\bs \{0\}\to \P({\cal H}), \ \ v\mapsto [v]$$
and we equip $\P({\cal H})$ with the quotient topology of  $p$. The
projective space 
 $\P({\cal H})$ becomes a Hilbert manifold as follows. 
For $v\in {\cal H}$ a unit vector we define an open neighborhood of $[v]$
in $\P({\cal H})$ by 

$$U_{[v]}=\{ [w]\in \P({\cal H})\: \la v,w\ra \neq 0\}.$$
Then the mapping 

$$\psi_{[v]}\: U_{[v]}\to v^\bot, \ \ [w]\mapsto {w\over \la w, v\ra }- v$$
is a homeomorphism with inverse 
$$ \psi_{[v]}^{-1}\: v^\bot\to U_{[v]}, \ \ u\mapsto [v+u].$$
A Riemannian metric on ${\cal H}\bs \{0\}$ is defined by 

$$g_v(u,w)={\Re \la u, w\ra\over \|v\|^2} \qquad (v\in {\cal H}\bs \{0\})
(u,w\in T_u ({\cal H}\bs \{0\})={\cal H}).$$

This metric is invariant under multiplication by $\C^*$ hence pushes 
down to a metric on $\P({\cal H})$. One obtains the familiar 
Fubiny-Study metric on $\P({\cal H})$ which we shall denote by 
$g^{FS}$. Explicitly it is given by 

$$g_{[v]}^{FS}\left( {d\over dt}\Big|_{t=0} [v+tu],  
{d\over dt}\Big|_{t=0} [v+tu]\right)={\|v\|^2 \la u, u\ra -|\la u,v\ra|^2
\over \|v\|^4}.\leqno(3.1)$$
We note that $(\P({\cal H}),g^{FS})$ is a metrically and geodesically 
complete space. 

\ssk
Associated to $\pi\in \hat G_s -\{ \1\}$ we consider the mapping 

$$F_\pi\: \Xi\to \P({\cal H}), \ \ xK_\C\mapsto [\pi(x)v_0].$$
This mapping is holomorphic as the orbit mapping $\Xi\to {\cal H}, 
\ \ xK_\C\mapsto \pi(x)v_0$ is holomorphic. 
Then we have the relation 
$$ g^\pi=2F_\pi^* g^{FS}\leqno(3.2)$$
which is a consequence of formula (3.1) and 
the standard computation: 
$$\eqalign{(\partial_Z\oline\partial_Z &
\log (s_\pi))(xK_\C) =\partial_Z \Big(2{\la\pi(x)v_0, d\pi(Z)\pi(x)v_0\ra\over 
\|\pi(x)v_0\|^2}\Big) \cr & =2{\la \pi(x)v_0, \pi(x)v_0\ra \la d\pi(Z)\pi(x)v_0, 
d\pi(Z)\pi(x)v_0\ra -|\la d\pi(Z)\pi(x)v_0,\pi(x)v_0\ra|^2\over\|\pi(x)v_0\|^4} 
.\cr}$$

\Proposition 3.1. Let $\pi\in \hat G_s-\{\1\}$. Then the map 
$F_\pi\: \Xi\to \P({\cal H}), \ xK_\C\mapsto [\pi(x)v_0]$ is proper. To be more precise, for every 
chart $\psi_{[v]}\: U_{[v]}\to v^\bot$ and every bounded closed ball 
$B\subeq v^\bot$ the subset $F_\pi^{-1}(\psi_{[v]}^{-1}(B))$
is compact in $\Xi$. 

\Proof. Fix a unit vector $v\in {\cal H}$ and a closed ball of finite radius
$B\subeq v^\bot$. 
Then for all $xK_\C\in F_\pi^{-1}(\psi_{[v]}^{-1}(B))$ we have 
$$\psi_{[v]}(F_\pi(xK_\C))={\pi(x)v_0\over \la \pi(x)v_0,v\ra} - v.\leqno(3.3)$$
To obtain  a contradiction let us assume that there exists a sequence 
$x_nK_\C$ in $F_\pi^{-1}(\psi_{[v]}^{-1}(B))$ such that  
$x_nK_\C\to \infty$. We can write $x_n=g_n\exp(iX_n)$ for $g_n\in G$ 
and $X_n\in \Omega$. As $\Omega$ is relatively 
compact, we may assume that $X_n\to X\in \oline \Omega$. 
{} From $x_nK_\C \to \infty$ we either deduce  $X\in \partial\Omega$ or 
$g_n\to \infty$. If $X\in \partial\Omega$, then 
Theorem 2.4 together with (3.3) shows that
$\psi_{[v]}(F_\pi(x_nK_\C))$ becomes unbounded. 
Hence $X_n\to X\in \Omega$ and so $g_n\to \infty$. Thus $\pi(\exp(iX_n))v_0$ lies in 
a compact subset $Q\subeq {\cal H}$. As $g_n\to \infty$, the 
Howe-Moore theorem on vanishing of matrix coefficients 
implies 
$$\lim_{n\to \infty} |\la \pi(x_n)v_0, v\ra|\leq
\lim_{n\to \infty}\sup_{w\in Q} |\la \pi(g_n^{-1})v, w\ra |=0.$$
In view of (3.3), this  shows that $\psi_{[v]}(F_\pi(x_nK_\C))$
becomes unbounded.  Again a contradiction, completing the proof. \qed

\Corollary 3.2.  Let $\pi\in \hat G_s-\{\1\}$. Then $\Im F_\pi$ 
is a closed submanifold of $\P({\cal H})$. In particular,
the $G$-invariant Riemannian metric $g^\pi$ associated to $\log (s_\pi)$ is complete.

\Proof. It follows from Proposition 3.1 that $\im F_\pi$ is a closed submanifold 
of the complete Riemannian manifold $(\P({\cal H}), g^{FS})$.  Now the 
assertion follows from (3.2). \qed

\msk Our discussion in this section leads to a formulation in complex-analytic terms of the non-trivial 
$K$-spherical unitary dual of $G$. We hope to expand on this result in a 
subsequent publication.

\Theorem 3.3. The map $\pi\to\omega_\pi = {i\over 2}\partial\oline\partial\log 
(s_\pi)$ identifies $\hat G_s -\{1\}$ with positive K\"ahler forms on $\Xi$ 
whose associated Riemannian metric is complete.

\Proof. The geometric properties of $\omega_\pi$ have been verified 
in Proposition 2.1 and Corollary 3.2. It remains 
to show that the assignment $\pi\mapsto\omega_\pi$ 
is injective. Let $\pi, \pi'\in \hat G_s-\{1\}$ such that 
$\omega_\pi=\omega_{\pi'}$. We recall that $\Xi$ 
is contractible (cf.\ Remark 4.5(b) below). Using Theorem 2.4 we showed in Proposition 2.6
that $\Xi$ is Stein.  Thus the potentials
$\log s_\pi$ and $\log s_{\pi'}$ for 
$\omega_\pi$ and $\omega_{\pi'}$ differ by the real part of a holomorphic 
function: 

$$\log s_\pi -\log s_{\pi'} =f\leqno(3.4)$$
where $f=\Re F$ for some $F\in {\cal O}(\Xi)$. Notice that $f$ is 
$G$-invariant as the left hand side of (3.4) is. In the sequel we 
identify $T(\Omega)$ with its image in $\Xi$. 
Set $H=F\res_{T(\Omega)}$ and $h=\Re H$. As $h$ is $A$-invariant 
pluriharmonic function we obtain: 

$$H(\exp(Z))=c+ i\lambda(Z) \qquad (Z\in \a+i\Omega)\ , $$ 
where $\lambda\in \a^*$ and $c\in \C$ is a constant.  
As $f$ is $G$-invariant, it follows in addition that 
$h$ is ${\cal W}$-invariant. Thus $\lambda=0$ and so $H=c$ is constant. 
It follows that $F=c$ is constant by the $G$-invariance of $f$. 

\par We conclude from (3.4) that 
$s_{\pi}=C s_{\pi'}$
for $C=e^{\Re c}>0$. Actually we have $C=1$ as $s_\pi(K_\C)=s_{\pi'}(K_\C)=1$. 
Thus 
$$s_{\pi}=s_{\pi'}.\leqno(3.5)$$
\par {}From (3.5) one can show that $\pi=\pi'$ as follows. 
For $\pi\in \hat G_s$ denote by 
$$\phi_\pi(gK)=\la \pi(g)v_0, v_0\ra\qquad  (g\in G)$$
the associated spherical function. This function extends holomorphically 
to $\Xi$ and its restriction to $A$ extends holomorphically 
to $T(2\Omega)$ (cf.\ [KS04, \S 4]). Denote these holomorphic extensions also 
by $ \phi_\pi$. The relation with $s_\pi$ is 
$$(\forall a\in \exp(i\Omega))\qquad \phi_\pi(a^2)=s_\pi(a).$$
Thus (3.5) implies $\phi_\pi=\phi_{\pi'}$ and so $\pi=\pi'$.
\qed 

\bobheadline{\S 4  Geometric analysis on $\Xi$}

In this section we will introduce a canonical $G$-invariant K\"ahler structure 
on $\Xi$. Lie theoretically, the K\"ahler structure is determined by the choice of 
Cartan-Killing form on $\g$, the Levi-Civita connection on G/K and their natural extensions to $\g_\Bbb C$ . It will be easy to see that they agree with the geometric description given in the introduction. We begin the 
section with a useful parametrization of $\Xi$.

\Proposition 4.1.  The map $$\Phi\:  G\times \Omega\to \Xi, \ \
(g,X)\mapsto g\exp(iX)K_\C$$ is a continuous surjection.  Moreover, 
$\Phi(g,X)=\Phi(g',X')$ for $g,g'\in G$ and $X,X'\in \Omega$ if and only if
there exists $k\in Z_K(X)$ and $w\in N_K(\a)$ such that $$g=g'wk\qquad
\hbox{and}\qquad X=\Ad(w^{-1}) X'.$$

\Proof.  Essentially this is Proposition 4 in [AG90].\qed

Write $\Omega'$ for the regular elements in $\Omega$, and $\Omega^+$ for
a connected component of $\Omega'$.  As usual we set $M=Z_K(\a)$ and
$\oline N=\theta(N)$.  Then we have the following refinement of
Proposition 4.1.

\Corollary 4.2.  \item{(i)} The map
$$\Phi'\:  G/M\times \Omega^+\to \Xi,\ \ (g,X)\mapsto g\exp(iX)K_\C$$ is an
analytic diffeomorphism onto its open and dense image $\Xi'$.  \item{(ii)}
The map $$\Phi''\:  \oline N A N\times \Omega^+, \ \ (g,X)\mapsto
g\exp(iX)K_\C$$ is an analytic diffeomorphism onto its open and dense image
$\Xi''$.

\Proof.  (i) Later, in the proof of Proposition 4.6, we will show that $\Phi'$ 
has everywhere regular differential.  (i) will then be a direct consequence of 
Proposition 4.1.
\par\nin (ii) In view of the Bruhat decomposition $$G=\coprod_{w\in {\cal W}}
\oline N w MA N,$$ with $\oline N MAN$ the open dense cell, (ii) follows 
immediately from (i).\qed

\subheadline{The canonical K\"ahler structure}

As before, for $X\in\g_\C$ we have defined a vector field $\tilde X$ on 
$G_\C/K_\C$ by

$$(\tilde X f)(x)={d\over dt}\Big|_{t=0} f(\exp(tX)x). $$ 
\par Let $p\:  G_\C \to G_\C /K_\C$ be the natural projection and denote by 
$p_\#\: \g_\C \to \p_\C $ the linear projection along $\k_\C$.  Let
$\rm{B}(\cdot,\cdot)$ be the Cartan-Killing form on $\g$ as well as its complex 
linear extension to
$\g_\C$.  With respect to the real form $\g$ of $\g_\C$ we have the complex
conjugation, $X\mapsto\oline X$.

\par Using Proposition 4.1 we write $z\in \Xi$ as $z=gaK_\C$,
for some $g\in G$ and $a=\exp(iX)$, $X\in \Omega$. Let $T_z\Xi$ denote the real 
tangent space to $\Xi$ at $z$ and $\tilde Y_z, \tilde Z_z\in T_z\Xi$. We claim 
that 

$$ (\tilde Y_z, \tilde Z_z)\mapsto \Re
\rm{B}(p_\#(\Ad(ga)^{-1}Y), p_\#(\oline{\Ad(ga)^{-1}Z}))$$ gives a well-defined
Hermitian bilinear form on $T_z\Xi$, which we shall denote 
$$ h_z\:  T_z\Xi \times T_z\Xi\to \R.$$  In fact, $\tilde Y_z=0$ means that 
$Y\in
\Ad(ga)^{-1}\k_\C$, so $p_\#(\Ad(ga)^{-1}Y)=0$.  To see that
the definition does not depend on the choice of $g$ and $a$ which represent
$z$, suppose that $z=g'a'K_\C$.  Then Proposition 4.1 gives that
$g'a'=gak$ for some $k\in K$.  Then

$$\eqalign{\rm{B}(p_\#(\Ad(g'a')^{-1}Y),&\ p_\#(\oline{\Ad(g'a')^{-1}Z})) =
\rm{B}(p_\#(\Ad(gak)^{-1}Y), p_\#(\oline{\Ad(gak)^{-1}Z}))\cr &
=\rm{B}(p_\#(\Ad(k)^{-1}\Ad(ga)^{-1}Y), p_\#(\oline{ 
\Ad(k)^{-1}\Ad(ga)^{-1}Z}))\cr
& =\rm{B}(\Ad(k)^{-1}p_\#(\Ad(ga)^{-1}Y), \Ad(k)^{-1}
p_\#(\oline{\Ad(ga)^{-1}Z}))\cr & =\rm{B}(p_\#(\Ad(ga)^{-1}Y),
p_\#(\oline{\Ad(ga)^{-1}Z})), \cr}$$ where we used the fact that $\Ad(k)^{-1}$
commutes with $p_\#$ and the complex conjugation, and that $\rm{B}$ is invariant
under $\Ad(k)^{-1}$.

Taking real and imaginary parts we get
$$h_z=g_z + i \omega_z,\leqno(4.1)$$ with $g_z=\Re h_z$ symmetric, and 
$\omega_z=\Im h_z$ skew-symmetric.

\ssk Recall the open dense subdomain $\Xi''$ from Corollary 4.2(ii).  Write
the elements $z\in \Xi''$ as $z=gaK_\C$ with $g\in \oline N A N$ and $a\in
\exp(i\Omega^+)$.  For every $X\in\g_\C$ we define an auxilliary vector
field $\hat X$ on $\Xi''$ by

$$(\hat X f)(z)={d\over dt}\Big|_{t=0} f(ga\exp(tX)K_\C).$$ That
$\hat X$ is well-defined follows from the uniqueness of the parametrization.  
Also
observe that

$$T_z\Xi=\{ \hat X_z \: X\in \p_\C\} \qquad (z\in \Xi'').$$

\Lemma 4.3.  For all $z\in \Xi''$ and $X, Y\in \g_\C$ we have

$$h_z(\hat X_z, \hat Y_z)=\rm{B}(p_\#(X), p_\#(\oline Y)).$$

\Proof.  This is immediate from the definition of $(h_z(\cdot,
\cdot))_{z\in\Xi}$ and the auxilliary vector fields $\hat X$.\qed

For every $h\in G$ we let $\lambda_h$ denote left translation by $h$ on $G_\C/ 
K_\C$, $$\lambda_h\:  G_\C/ K_\C\to G_\C/ K_\C,\ \ x\mapsto hx.$$

\Theorem 4.4.  The Hermitian metric $h_z(\cdot,\cdot)$ on $\Xi$ is a 
$G$-invariant K\"ahler metric.

\Proof. Clearly this structure is smooth so from the above $h_z(\cdot,\cdot)$ is 
an Hermitian metric. To show $\big(h_z(\cdot,\cdot)\big)_{z\in\Xi}$ is 
$G$-invariant we have to show that

$$h_{hz}(d\lambda_h (\tilde Y_z), d\lambda_h (\tilde Z_z))=h_z(\tilde Y_z , 
\tilde
Z_z)$$ holds for all $z\in \Xi$ and $h\in G$.  Since $d\lambda_h (\tilde
Y_z)=\tilde {\Ad(h)Y}_{hz}$ and $d\lambda_h (\tilde Z_z)=\tilde {\Ad(h)Z}_{hz}$,
we have

$$\eqalign{h_{hz}(d\lambda_h (\tilde Y_z), d\lambda_h (\tilde Z_z)) 
&=h_{hz}(\tilde
{\Ad(h)Y}_{hz}, \tilde {\Ad(h)Z}_{hz})\cr &=\rm{B}(p_\#(\Ad(hga)^{-1}\Ad(h)Y),
p_\#(\oline{\Ad(hga)^{-1}\Ad(h)Z}))\cr &=\rm{B}(p_\#(\Ad(ga)^{-1}Y),
p_\#(\oline{\Ad(ga)^{-1}Z}))\cr &=h_z(\tilde Y_z, \tilde Z_z).\cr }$$

\par\nin Next we show the structure $\big(g_z(\cdot,\cdot)\big)_{z\in\Xi}$ is
Riemannian. By $G$-invariance it is enough to show for every $a\in 
\exp(i\Omega)$ that

$$g_a\:  T_a\Xi\times T_a\Xi\to \R$$ is positive definite, where, by abuse
of notation, we write $a$ for $aK_\C\in \Xi$.  For that, first we observe that

$$T_a\Xi=\{ \tilde Y_a\:  Y\in \p_\C\}.$$ Indeed, every $Y\in\p_\C$ can be
written as $$Y=Y_0+\sum_{\alpha\in \Sigma^+} (Y_\alpha -\theta Y_\alpha)$$
for $Y_0\in\a_\C$ and $Y_\alpha\in \g_\C^\alpha$.  So

$$\eqalign{\Ad(a)^{-1}Y &= Y_0+ \sum_{\alpha\in \Sigma^+} ( a^{-\alpha}
 Y_\alpha -a^{\alpha} \theta Y_\alpha)\cr &=Y_0+\sum_{\alpha\in \Sigma^+}
{a^\alpha+a^{-\alpha}\over 2} ( Y_\alpha - \theta Y_\alpha) +
\underbrace{\sum_{\alpha\in \Sigma^+} {a^{-\alpha}-a^\alpha\over 2} (
Y_\alpha + \theta Y_\alpha)}_{\in \k_\C}.\cr}$$ But then
$$p_\#(\Ad(a)^{-1}Y) = Y_0+\sum_{\alpha\in \Sigma^+} {a^\alpha+a^{-\alpha}\over
2} ( Y_\alpha - \theta Y_\alpha).  $$ Using the definition of $\Omega$ we get
for any $a=\exp(iX)$, $X\in \Omega$, $$ {a^\alpha+ a^{-\alpha}\over
2}=\cos\alpha(X)>0.$$ From this it follows that $T_a\Xi=\{ \tilde Y_a\:  Y\in 
\p_\C\}.$ In addition, for $Y\neq0$ we get an expression for $g_a(\tilde Y_a, 
\tilde Y_a)$ showing that $g_a(\cdot,\cdot)$ is positive definite, namely 
$$\eqalign{g_a(\tilde Y_a, \tilde Y_a)&=\Re\rm{B}(Y_0, \oline Y_0)+ 
\sum_{\alpha\in\Sigma^+} \Big({a^\alpha+ a^{-\alpha}\over 2}\Big)^2 \Re\rm{B}( 
Y_\alpha-\theta Y_\alpha, \oline Y_\alpha -\theta \oline Y_\alpha)\cr
&=\Re\rm{B}(Y_0, \oline Y_0)+ 2\sum_{\alpha\in\Sigma^+} \Big({a^\alpha+
a^{-\alpha}\over 2}\Big)^2 \Re \rm{B}( Y_\alpha, -\theta \oline Y_\alpha)>0
. \cr}$$

\ssk\nin Finally, to see that $\big(h_z(\cdot,\cdot) \big)_{z\in \Xi}$ is
K\"ahler we have to show that $\big(\omega_z(\cdot,\cdot)
\big)_{z\in\Xi}$ is closed.  Recall the open dense subdomain $\Xi''$ and the
auxilliary vectorfields $\hat X$ on $\Xi''$ for $X\in \g_\C$.  From Lemma 4.3
we have $$\omega_z(\hat X_z ,\hat Y_z)=\Im \rm{B}(p_\#(X), p_\#(\oline Y)).  
\leqno(4.2)$$ By the density of $\Xi''$ in $\Xi$, it suffices to show
that $d\omega_z(\hat X_z, \hat Y_z, \hat Z_z) =0$ for all $z\in \Xi''$ and $X,
Y, Z\in \p_\C$.  Now

$$\eqalign{d\omega(\hat X, \hat Y,\hat Z)&= \hat X\omega(\hat Y, \hat Z)
-\hat Y\omega(\hat X, \hat Z) +\hat Z\omega(\hat X, \hat Y)\cr & - \omega
([\hat X, \hat Y], \hat Z) + \omega ([\hat X, \hat Z], \hat Y) -\omega ([\hat
Y, \hat X], \hat Z).\cr}$$  On $\Xi''$ the first three terms on the right hand 
side
are zero by (4.2), while the last three terms are zero by (4.2) and
$[\p_\C,\p_\C]\subeq \k_\C$.  \qed

\Remark 4.5.  (a) One can see from the previous proof that $\big(g_z(\cdot, \cdot)\big)_{z\in \Xi}$extends the $G$-invariant metric on $G/K$.  \par\nin (b) The Riemannian structure $\big(g_z(\cdot, \cdot)\big)_{z\in \Xi}$
is the K\"ahler metric coming from the adapted complex structure. The adapted complex structure associated to a complete Riemannian manifold $M$ can be characterized as the unique complex structure on a neighborhood of the zero section of the tangent bundle of $M$ satisfying: for any geodesic $\gamma\:\R\to M, \ \ t\mapsto \gamma(t)$ the map

$$\phi_\gamma\:  \C \to TM, \ \ t+is\mapsto (\gamma(t), s\gamma'(t))$$ is holomorphic. 
\par For $M=G/K$ we have the isomorphism

$$T(G/K)=G\times_K\p$$ and the analytic $G$-equivariant map

$$\Phi\:  T(G/K)\to G_\C/ K_\C, \ \ [g,X]\mapsto g\exp(iX)K_\C.$$ If
$\gamma\:  \R\to G/K$ is a geodesic then $\gamma(t)=g\exp(tX)K$ for some
$g\in G$ and $X\in \p$.  So
$$\Phi(\phi_\gamma(s,t))=\Phi([\gamma(t), sX])=g\exp((t+is)X)K_\C.$$ In
particular, for all geodesics $\gamma$ the maps $\Phi\circ \phi_\gamma$ are
holomorphic.

\par Define $W=\Ad(K)(\Omega)$ and note that $W\subeq \p$ is a convex open
set (easy consequence from the Kostant convexity theorem).  Recall from
[AG90, Prop.\ 4] that the restriction of $\Phi$ to $G\times_K W$

$$G\times_K W\to \Xi, \ \ [g,X]\mapsto g\exp(iX)K_\C$$ is an analytic
diffeomorphism.  From our calculation above it now follows that $\Xi$ carries
the adapted complex structure, see also [BHH03].

\par\nin (c) The metric is not complete, since for every $X\in\oline \Omega$,
the curve $\gamma(t)=\exp(itX)$, $t\in [0,1[$ has finite length.  \qed

\subheadline{The Riemannian measure}

For the remainder of this section we consider $\Xi$ as a Riemannian manifold 
equipped with the Riemannian structure $(g_z(\cdot, \cdot))_{z\in\Xi}$.  We 
write $\mu_\Xi$ for the Riemannian measure on $\Xi$. The next proposition gives 
a simple description of the Riemannian measure in terms of coordinates. For this 
purpose, for $\alpha\in\Sigma$ set $m_\alpha=\dim \g^\alpha$. We remind the 
reader of the notation $\Omega^+\subeq \Omega$ for the regular chamber.

\Proposition 4.6. For $f\in L^1(\Xi)$ 

$$\int_\Xi f(z)\ d\mu_\Xi(z)=c\int_{G/M} \int_{\Omega^+} f(g\exp(iX)K_\C) \
\prod_{\alpha\in\Sigma^+} |\sin 2\alpha(X)|^{m_\alpha} \ dX \ d\mu_{G/M}(gM)
,$$
where $c>0$ depends only on the normalization of the Haar measure $\mu_{G/M}$ on 
$G/M$. 

\Proof.  Recall that the parametrization

$$\Phi'\:  G/M\times \Omega^+ \to \Xi, \ \ (g,X)\mapsto g\exp(iX)K_\C$$ is a 
smooth injective map.  We show that $\Phi'$ is everywhere
regular and we compute the Jacobian. For $gM\in G/M$ and $X\in \Omega^+$
we write $z=g\exp(iX)K_\C$. We use the following
identifications of tangent spaces: $$T_{(gM,X)}(G/M\times \Omega^+)\simeq
\g/\m\times \a\simeq \m^\bot \times \a$$ and $$T_z\Xi\simeq\p_\C,  \quad
\hat Y_z\leftrightarrow Y \qquad (Y\in\p_\C).$$ 
With these and Lemma 4.3 we
describe an inner product on $\m^\bot\times \a$ (independent of $g$ and $X$)
such that the map

$$T_{(gM,X)}(G/M\times \Omega^+)\to T_z\Xi=\p_\C$$
$$(Y,Z)\mapsto \cases { Y+iZ\ & for $Y\in\a$ \cr Y-\theta Y& for $Z=0$ and
$Y\in \g^\alpha$, $\alpha\in \Sigma^+$\cr iY-i\theta Y& for $Z=0$ and $Y\in
\g^\alpha$, $\alpha\in \Sigma^-$\cr}\leqno(4.3)$$ becomes an isometric
isomorphism of tangent spaces.  With respect to this basis and identifications 
we compute the Jacobian $|\det d\Phi(gM,X)|$ of the differential
$d\Phi(gM,X)$.  By $G$-equivariance we have $$|\det d\Phi(gM,X)|=|\det
d\Phi(\1,X)|,$$ so we may assume that $g=\1$.  Fix $X\in \Omega^+$, set
$a=\exp(iX)$ and $z=aK_\C$.  Then for $Y\in\m^\bot$ and
$Z\in \a$ we have $$\eqalign{d\Phi(\1, X)(Y,Z)& = {d\over dt}\Big|_{t=0}
\exp(tY)\exp(X)K_\C + {d\over dt}\Big|_{t=0} \exp(iX)\exp(itZ)K_\C\cr &
=\hat{\Ad(a)^{-1}Y}_z + \hat {iZ}_z \leftrightarrow \Ad(a)^{-1}Y +iZ.
\cr}\leqno(4.4)$$

As $Y$ is in $\m^\bot$ it can be written $Y=Y_0 +\sum_{\alpha\in\Sigma} 
Y_\alpha$ with
$Y_0\in \a$ and $Y_\alpha\in\g^\alpha$.  For any $W\in\p_\C$

$$\eqalign{g_z(d\Phi(\1, X)(Y,Z), \hat W_z)&=\Re \rm{B} \big(p_\#(\Ad(a)^{-1}Y +
iZ), W\big)\cr &=\Re \rm{B} \big(p_\#(\Ad(a)^{-1}Y + iZ), W\big)\cr &=\Re
\rm{B}(iZ,W) + \Re \rm{B}(Y_0,W) + \sum_{\alpha\in\Sigma} \Re \rm{B}
\big(p_\#(a^\alpha Y_\alpha), W\big).\cr}\leqno(4.5)$$ Putting together the
identifications of tangent spaces, (4.3), (4.4) and (4.5) one obtains $$|\det 
d\Phi(\1,X)|=c \prod_{\alpha\in\Sigma^+} \big|\det
d\Phi(\1,X) \res_{(\g^\alpha\oplus\g^{-\alpha})\times\{0\}}\big|
.\leqno(4.6)$$

With $Y_\alpha\in \g^\alpha$ a unit vector (4.3)-(4.5) give

$$\eqalign{\big|\det d\Phi(\1,X)&
\res_{(\g^\alpha\oplus\g^{-\alpha})\times\{0\}}\big|=\cr &=\det \Big|\pmatrix
{\Re \rm{B} \big (p_\#(a^{-\alpha} Y_\alpha), Y_\alpha -\theta Y_\alpha\big) &
\Re \rm{B} \big (p_\#(a^{-\alpha} Y_\alpha), iY_\alpha -i\theta Y_\alpha\big)\cr
\Re \rm{B} \big (p_\#(a^\alpha \theta Y_\alpha), Y_\alpha -\theta Y_\alpha\big)
& \Re \rm{B} \big (p_\#(a^\alpha \theta Y_\alpha), iY_\alpha -i\theta
Y_\alpha\big)\cr}\Big|^{m_\alpha}\cr &=\det \Big|\pmatrix{\Re a^{-\alpha} &
\Re ia^{-\alpha}\cr \Re a^\alpha & \Re ia^\alpha\cr}\Big|^{m_\alpha} =\det
\Big|\pmatrix{\cos \alpha(X) & \sin \alpha(X)\cr \cos\alpha(X) & -\sin
\alpha(X)\cr}\Big|^{m_\alpha}\cr &=|\sin 2\alpha(X)|^{m_\alpha} .  \cr}$$
So from (4.6) we get

$$|\det d\Phi(\1,X)|=c\prod_{\alpha\in \Sigma^+} |\sin
2\alpha(X)|^{m_\alpha}\neq 0,$$ concluding the proof of the
proposition.\qed

\Remark 4.7.  From the formula in Proposition 4.6 one sees that the Riemannian 
measure $\mu_\Xi$ equals the restriction of a Haar measure $\mu_{G_\C/ K_\C}$ to 
$\Xi$. Similarly, one sees that the Riemannian measure $\mu_\Xi$ vanishes at 
infinity, i.e. the density vanishes at the boundary $\partial\Xi$.  This, and 
the precise order of the decay, will be used later.\qed

Our next goal is to identify the Laplace-Beltrami operator on $\Xi$ 
corresponding to this Riemannian structure.  For that we have to recall some 
facts about invariant differential operators
on $G_\C/ K_\C$.

\subheadline{Invariant differential operators on $G_\C/ K_\C$}

We begin with some standard facts about the complexification of complex Lie 
algebras. Let $\g^\C_\R$ denote $\g_\C$ considered as a real Lie algebra. 
$\g^\C_\R$ has two natural (almost) complex structures corresponding to 
multiplication by $i$ and by $-i$. With respect to the first $\g^\C_\R$ is 
canonically isomorphic to $\g_\C$. Denote by $\oline\g_\C$the complex Lie 
algebra $\g_\C$ equipped with the conjugate complex structure. Then $\g^\C_\R$ 
with the second complex structure is canonically isomorphic to $\oline\g_\C$. 
Let$\g_\C^\C$ denote the complexification of the real Lie algebra $\g^\C_\R$.  
As a complex vector space $\g_\C^\C$ has a complex structure, say$J_{\C}$. 
Relative to these structures the map

$$(\g_\C^\C, J_{\C})\to \g_\C\oplus\oline\g_\C, \ \ X+J_{\C}Y\mapsto (X+iY,
X-iY)\leqno(4.7)$$ is an isomorphism of complex Lie algebras.

Let $\Diff_L(G_\C)$ denote the left $G_\C$-invariant differential
operators on $G_\C$. Any $X\in \g_\C$ gives a left invariant differential 
operator $L_X$ by

$$(L_Xf)(x)={d\over dt}\Big|_{t=0} f(g\exp(tX)),$$ for $x\in G_\C$ and
$f$ a function on $G_\C$ differentiable at $x$.

Then the map $$\g_\C^\C\to \Diff_L(G_\C), \ \ X+J_{\C}Y\mapsto L_X + iL_Y$$
extends to an algebra isomorphism $${\cal U}(\g_\C^\C)\simeq \Diff_L(G_\C).$$
As $K_\C$ is reductive in $G_\C$, a left invariant differential
operator on $G_\C$, viewed as an element in ${\cal U}(\g_\C^\C)$, projects to
an operator on $G_\C/ K_\C$ if and only if it is $\Ad(K_\C)$-invariant. In other 
words,

$${{\cal U}(\g_\C^\C)^{K_\C}\over {\cal
U}(\g_\C^\C)^{K_\C}\cap {\cal U}(\g_\C^\C)\k_\C^\C}\simeq \Diff_L(G_\C/ 
K_\C).\leqno(4.8)$$

Let $X_1,\ldots, X_n$ be an orthonormal basis of $\p$ and set
$Y_j=iX_j\in i\p$. With $$Z_j={1\over 2}(X_j -J_{\C}
Y_j)\qquad \hbox{and}\qquad \oline Z_j={1\over 2} (X_j+ J_{\C}Y_j)$$ and using 
the identification (4.7) we set, as before, $$L_{Z_j}={1\over 
2}(L_{X_j}-iL_{Y_j}) \qquad \hbox{and} \qquad L_{\oline Z_j}={1\over 2}(L_{X_j} 
+i L_{Y_j}).$$ 
Consider the differential operators on $G_\C$ defined by

$$\Delta=\sum_{j=1}^n L_{X_j}^2, $$ $$\Delta_\C=\sum_{j=1}^n L_{X_j}^2
+\sum_{j=1}^n L_{Y_j}^2, $$ $$\square=\sum_{j=1}^n L_{Z_j}^2, $$ and
$$\oline\square=\sum_{j=1}^n L_{\oline Z_j}^2.$$

\Proposition 4.9.  Each of the differential operators $\Delta$, $\Delta_\C$,
$\square$ and $\oline \square$ induces a differential operator on $G_\C/
K_\C$.

\Proof.  We consider first the operator $\Delta$.  Since $X_1,\ldots,
X_n$ is an orthonormal basis of $\p$, $\sum_{j=1}^n
X_j^2\in {\cal U}(\g)$ is $\Ad(K)$-invariant.  In particular, $\sum_{j=1}^n
X_j^2\in {\cal U}(\g_\C)^K$ and, since $K_\C$ is reductive, $$\sum_{j=1}^n 
X_j^2\in {\cal
U}(\g_\C)^{K_\C}.$$  In view of
(4.8), it follows that $\Delta$ projects to a differential operator on
$G_\C/ K_\C$.  One then shows mutatis mutandis that $\Delta_\C$ factors
to $G_\C/ K_\C$.

\par For the operators $\square$ and $\oline \square$
we use the identification ${\cal U}(\g_\C^\C)\simeq {\cal U}(\g_\C) \times
{\cal U}(\oline\g_\C)$ induced by (4.7).  Then $\square $ and $\oline
\square$ become identified with the operators

$$\square \leftrightarrow \big(0, \sum_{j=1}^n X_j^2\big)\qquad
\hbox{and}\qquad \oline\square \leftrightarrow\big( \sum_{j=1}^n X_j^2,
0\big).  \leqno(4.9)$$ As observed above it follows that $\square$ and $\oline 
\square$ are represented by
$\Ad(K_\C)$-invariant elements in ${\cal U}(\g_\C^\C)$, hence induce 
differential operators on $G_\C/ K_\C$.\qed

By abuse of notation we use $\Delta$, $\Delta_\C$, $\square$ and $\oline 
\square$
to denote either differential operators on $G_\C/ K_\C$ or on $G_\C$. Their use 
should be clear from the context.

\Proposition 4.10.  For $\Delta$, $\Delta_\C$, $\square$ and $\oline \square$
considered as elements in $\Diff_L(G_\C/ K_\C)$ the following assertions
hold:

\item{(i)} $[\square ,\oline \square]=0$;  \item{(ii)} $[\Delta_\C, \square+
\oline \square]=0$;  \item{(iii)} $\Delta={1\over 2}\Delta_\C +(\square
+\oline \square)$.

\Proof.  (i) This is immediate from (4.9).  \par\nin (ii) In view of (4.8),
the assertion in (ii) can be phrased as

$$[\sum_{j=1}^n X_j^2 +\sum_{j=1}^n Y_j^2, \sum_{j=1}^n Z_j^2 + \sum_{j=1}^n
\oline Z_j^2]\in {\cal U}(\g_\C^\C)\k_\C^\C.\leqno(4.10)$$ An easy
computation shows that $$\square +\oline \square={1\over 2}(\sum_{j=1}^n
L_{X_j}^2 - \sum_{j=1}^n L_{Y_j}^2)\leftrightarrow {1\over
2}\big(\sum_{j=1}^n X_j^2 - \sum_{j=1}^n Y_j^2\big).\leqno(4.11)$$ Substituting 
(4.11) into the left side of (4.10) we obtain
$$\eqalign{[\sum_{j=1}^n X_j^2 +\sum_{j=1}^n Y_j^2, \sum_{j=1}^n Z_j^2 +
\sum_{j=1}^n \oline Z_j^2]& ={1\over 2} [\sum_{j=1}^n X_j^2 +\sum_{j=1}^n
Y_j^2, \sum_{j=1}^n X_j^2 -\sum_{j=1}^n Y_j^2]\cr &= [\sum_{j=1}^n X_j^2,
\sum_{j=1}^n Y_j^2]=\sum_{j,k=1}^n[X_j^2, Y_k^2].  \cr} \leqno(4.12)$$ Since
$X_j, Y_k\in\p_\C$, $[X_j, Y_k]\in \k_\C$. An easy computation
shows that $[X_j^2, Y_k^2]\in {\cal U}(\g_\C) \k_\C$ and so (4.12) implies
(4.10), completing the proof of (ii).  \par\nin (iii) This is immediate from
(4.11) and the definition of $\Delta$ and $\Delta_\C$.  \qed

\Remark 4.11.  The operators $\square $, $\oline \square$ and $\Delta$ can be
related through the previously described holomorphic extension result. Indeed, 
for $x\in G/K$ let $f$ be a function defined in a connected neighborhood $U$ of 
$x$
which extends holomorphically to a connected complex neighborhood $U_\C\subeq
G_\C/ K_\C$ with $U_\C \cap G/K=U$.  Write $f^\sim$ for the holomorphic
extension of $f$ to $U_\C$.  Then $\square$, $\oline \square$ and $\Delta$
are related through $$(\Delta f)^\sim(z)=(\square f^\sim)(z) \qquad (\forall
z\in U_\C)$$ and $$\oline{(\Delta f)^\sim}(z)= (\oline \square\,
\oline{f^\sim})(z) \qquad (\forall z\in U_\C).\qeddis

\subheadline{The Laplace-Beltrami operator on $\Xi$}

We write $\Delta_\Xi$ for the restriction of $\Delta_\C$ to
$\Xi$.  The next result establishes the naturality of the previously defined 
Riemannian metric on $\Xi$.

\Theorem 4.12.  The operator $\Delta_\Xi$ is the Laplace-Beltrami operator
for the Riemannian metric $\big(g_z(\cdot, \cdot)\big)_{z\in \Xi}$ on $\Xi$.

\Proof.  We use the open dense subdomain $\Xi''$.  In particular, it follows
from Lemma 4.3 that for $z=gaK_\C\in \Xi''$ and $X,Y\in\p_\C$ we have $$g_z(\hat
X_z ,\hat Y_z)=\Re \rm{B}(X, \oline Y).$$ Thus with $X_1,\ldots, X_n$
and $Y_1, \ldots, Y_n$ as before, the Laplace-Beltrami operator on
$\Xi''$ is given by

$$\sum_{j=1}^n \hat X_j^2 +\sum_{j=1}^n\hat Y_j^2.$$ But this operator
evidently coincides with $\Delta_\Xi\res_{\Xi''} =\Delta_\C\res_{\Xi''}$.
Since $\Xi''$ is open dense, and $\Delta_\C$ has analytic coefficients, the 
theorem follows.\qed

\bobheadline{\S 5  Invariant Hilbert spaces of holomorphic functions on
$\Xi$}

The holomorphic extension result in Theorem 1.1 has as a consequence a 
realization of each representation in $\hat G_s$ as a Hilbert space of 
holomorphic functions on $\Xi$. In this section we draw some natural conclusions 
from this point of view. We shall use some facts about Hilbert spaces of 
holomorphic functions for which we refer the reader to [FT99] for a recent 
account.
 
For $M$ a second countable complex manifold, we write ${\cal O}(M)$ for the 
Fr\'echet space of holomorphic functions on $M$. A Hilbert space of holomorphic 
functions on $M$ is aHilbert space ${\cal H}\subeq {\cal O}(M)$ such that the
inclusion mapping ${\cal H}\into {\cal O}(M)$ is continuous. For an equivariant 
version suppose that $G$ is alocally compact group that acts by holomorphic 
automorphisms on $M$. One says that ${\cal H}$ is a {\it$G$-invariant Hilbert 
space of holomorphic functions on $M$} if the left regular representation of $G$ 
on ${\cal H}$ is a unitary representation of $G$.

If ${\cal H}\subeq {\cal O}(M)$ is a Hilbert space of holomorphic functions, 
then the point evaluation map, $f\mapsto f(z)$, is continuous for all $z\in M$. 
In the usual way one obtains a reproducing kernel for ${\cal H}$  $$K\:  M\times
M\to\C,$$  
holomorphic in the first variable and anti-holomorphic in the second variable. 
Unitarity of $G$ on ${\cal H}$ becomes $G$-invariance for $K$, $K(gz, 
gw)=K(z,w)$ for all $g\in G$ and $z,w\in M$.

\par As is well known, this use of the Riesz representation theorem allows one 
to establish an equivalence between $G$-invariant Hilbert spaces of holomorphic 
functions on $M$ and such kernels. Indeed, if we let $\oline M$ denote $M$ with 
the conjugate complex structure, then the kernel $K$ is in ${\cal O} (M\times 
\oline M)$, and, if ${\cal P}_G(M)\subeq {\cal O} (M\times \oline M)$ is the 
cone of $G$-invariant holomorphic positive definite kernels on $M$, then $K\in 
{\cal P}_G(M)$. Now ${\cal P}_G(M)$ can be shown to be a closed convex conuclear 
cone in the Fr\'echet space ${\cal O} (M\times \oline M)$. Let $\Ext({\cal
P}_G(M))$ be the cone of extremal rays in ${\cal P}_G(M)$. The elements in 
$\Ext({\cal P}_G(M))$ correspond to the irreducible $G$-invariant Hilbert spaces 
of holomorphic functions on $M$.Then reducibility questions concerning the 
action of $G$ on ${\cal H}$ are re-formulated into so-called Choquet theory.

\par Now we specialize this by taking for $M$ the domain $\Xi\subeq G_\C/ K_\C$ 
on
which $G$ acts by holomorphic automorphisms.

\par As before, let $\hat G_s$ denote the $K$-spherical unitary dual of $G$.  
Let$[\pi]\in \hat G_s$ and $(\pi, {\cal H}_\pi)$ a representative of it with
$v_0\in {\cal H}_\pi$ a normalized $K$-spherical vector.  Then Corollary
1.2 implies the existence of a continuous $G$-equivariant embedding
$$\iota_\pi\:  {\cal H}_\pi\into {\cal O}(\Xi), \ \ v\mapsto f_v; \
f_v(xK_\C)=\la \pi(x^{-1})v, v_0\ra .\leqno(5.1)$$ In this way $\iota_\pi({\cal
H}_\pi)$ equipped with the topology of ${\cal H}_\pi$ becomes a $G$-invariant
Hilbert space of holomorphic functions.  {} From (5.1) it follows that the
reproducing kernel $K_\pi$ of $\iota_\pi({\cal H}_\pi)$ is given by
$$K_\pi\:  \Xi\times\Xi\to \C, \ \ (xK_\C ,yK_\C)\mapsto \la \pi(\oline y)
v_0, \pi(\oline x)v_0\ra ,\leqno(5.2)$$ where $x\mapsto \oline x$ denotes
the complex conjugation on $G_\C$ with respect to the real form $G$. The 
re-formulation of $\hat G_s$ into kernels as provided by the general theory  
gives $$\hat G_s\leftrightarrow\Ext({\cal P}_G(\Xi))=\coprod_{[\pi]\in \hat G_s} 
\R^+K_\pi .$$

The decomposition of invariant Hilbert spaces of holomorphic 
functions on $\Xi$ into irreducibles is a direct consequence of results in 
[FT99] or [K99]. For this reason we simply sketch the argument. Independently Faraut [F03] made the same observation with the same argument in a lecture at MSRI, '01. 

\par Identifying $\hat G_s$ with a subset of $\a_\C^*$, we write
$(\pi_\lambda, {\cal H}_\lambda)$ for a representative of $\lambda\in \hat
G_s\subeq \a_\C^*$. Similarly we write $K_\lambda$ for
$K_{\pi_\lambda}$.  It is known that the Borel structure on $\hat G_s$ induced
from the hull-kernel topology on $\hat G$ is the same as the Borel structure on 
$\a_\C^*$ induced from the Euclidean topology. Also, the map $$\hat G_s\to 
\Ext({\cal P}_G(\Xi)), \ \
\lambda\mapsto K_\lambda$$ is continuous and hence constitutes an {\it
admissible parametrization} of $\Ext({\cal P}_G(\Xi))$ in the sense of
[K99].  The following result is an immediate consequence of the abstract 
results in [FT99] or [K99].

\Theorem 5.1.  Let ${\cal H}$ be a $G$-invariant Hilbert space of holomorphic
functions on $\Xi$.  Let $K(z,w)$ be the reproducing kernel of ${\cal H}$.
Then there exists a unique Borel measure $\mu$ on $\hat G_s$ such that
$$K(z,w)=\int_{\hat G_s} K_\lambda(z,w) \ d\mu(\lambda) \qquad (z,w\in \Xi)$$
with the right hand side converging absolutely on compact subsets of
$\Xi\times\Xi$.\qed

Theorem 5.1 implies a unitary $G$-equivalence of the left regular representation 
 of $G$ on ${\cal H}$ with a direct integral of representations from $\hat G_s$ 
relative to a spectral measure $\mu$

$$(L, {\cal H})\simeq \Big(\int_{\hat G_s}^\oplus \pi_\lambda \ d\mu(\lambda),
\int_{\hat G_s}^\oplus {\cal H}_\lambda \ d\mu(\lambda)\Big).
\leqno(5.3)$$

\par Shortly we shall give a criterion for a Borel
measure $\mu$ on $\hat G_s$ to be a spectral measure of a $G$-invariant
Hilbert space of holomorphic functions on $\Xi$. But first we present some 
examples of invariant Hilbert spaces on $\Xi$.

\Example 5.2.  (a) The $G$-invariance of the Riemannian measure
$\mu_\Xi$ allows one to define the Bergmann space of $\Xi$, 
$${\cal B}^2(\Xi)=\{ f\in {\cal O}(\Xi)\:  \|f\|^2 =\int_\Xi |f(z)|^2\
d\mu_\Xi(z)<\infty\}.$$ Certainly ${\cal B}^2(\Xi)$ is a $G$-invariant Hilbert 
space of holomorphic functions on $\Xi$.  We conjecture that 
that  ${\cal B}^2(\Xi) \neq \{0\}$.
\par As usual, more generally one can define for every $G$-invariant 
weight function $p\: \Xi\to \R^+$ a {\it weighted} Bergman space
${\cal B}^2(\Xi, p)$. Again, these weighted Bergman spaces
are $G$-invariant Hilbert spaces of holomorphic functions 
on $\Xi$. A characterization of the spectral measure follows from Faraut's Gutzmer Formula ([F02,03]). 

\par\nin (b) Using notation to be introduced in \S 7, we suppose that 
$\Xi=\Xi_0$ and that the distinguished boundary
of $\Xi$ (cf.\ [GK02a]) is a symmetric space $G/H$. For such spaces one can 
define a Hardy space${\cal H}^2(\Xi)$ on $\Xi$. Then ${\cal H}^2(\Xi)$ is a 
$G$-invariant Hilbert-space of holomorphic functions on $\Xi$.  This Hardy space 
has a boundary value map ${\cal H}^2(\Xi)\to L^2(G/H)$ whose image is a full 
multiplicity one subspace of the most continuous spectrum of $L^2(G/H)$. The 
spectral measure here is nothing other than the contribution to the most 
continuous spectrum of the Plancherel measure. For all this see [GK\'O03, 04].\qed

\par We shall present a necessary condition for a Borel measure to be a spectral 
measure for a $G$-invariant Hilbert space of holomorphic functions on $\Xi$. We 
start with a simple observation.

\Lemma 5.3.  Let $\lambda\in \hat G_s$.  Then for all $g\in G$ and
$a=\exp(iX)$, $X\in\Omega$  we have  \item{(i)} if
$v\in {\cal H}_\lambda$, then $$|f_v(gaK_\C)|\leq \|\pi_\lambda(a^{-1})v_0\|
\cdot \|v\|;$$ \item{(ii)} $K_\lambda(gaK_\C,
gaK_\C)=\|\pi_\lambda(a^{-1})v_0\|^2.$

\Proof.  (i) Using the unitarity of $\pi_\lambda$ and the Cauchy-Schwarz
inequality one gets

$$\eqalign{|f_v(gaK_\C)|&=|\la \pi_\lambda(a^{-1}g^{-1})v, v_0\ra |= |\la
\pi_\lambda(g^{-1})v, \pi_\lambda(a^{-1})v_0\ra| \cr &\leq
\|\pi_\lambda(g^{-1})v\|\cdot \|\pi_\lambda(a^{-1})v_0\|=
\|\pi_\lambda(a^{-1})v_0\| \cdot \|v\|.\cr}$$  
\par\nin (ii) From the definition of $K_\lambda$ (cf.\ (5.2)) we have

$$K_\lambda(gaK_\C, gaK_\C)=\la\pi_\lambda(ga^{-1})v_0,
\pi_\lambda(ga^{-1})v_0\ra=\la\pi_\lambda(a^{-1})v_0,
\pi_\lambda(a^{-1})v_0\ra=\|\pi_\lambda(a^{-1})v_0\|^2, $$ as was to be
shown.\qed

Let us define a norm $|\cdot|_\Omega$ on $\a^*$ by 
$$|\lambda|_\Omega\:=\sup_{X\in\Omega} |\lambda(X)| \qquad (\lambda\in \a^*).$$
Recall from [KS04, Sect.\ 4] that if $\phi_\lambda$ is the analytically 
continued spherical function on $G/K$ with parameter $\lambda$, and $a\in 
\exp(i\Omega)$,

$$\|\pi_\lambda(a^{-1})v_0\|^2=\phi_\lambda(a^2).\leqno(5.4)$$

\par Let $Q\subeq \Omega$ be a compact subset.  Then it follows from
(5.4), (1.1) and [KS04, Th.\ 4.2] that there exist constants $C_Q>0$ and 
$\epsilon_Q,\ 0\le\epsilon_Q < 1$ such
that

$$\|\pi_\lambda(a^{-1})v_0\|^2\leq C_Q e^{(1 -\epsilon_Q)2|\Im\lambda|_\Omega}\leqno 
(5.5)$$
holds for all $a\in \exp(iQ)$.

\Proposition 5.4.  Let $\mu$ be a Borel measure on $\hat G_s$ such that
$$(\forall 0\leq c<2)\qquad \int_{\hat G_s} e^{c|\Im \lambda|_\Omega} \ d\mu(\lambda)<\infty .$$ Then
$\mu$ is the spectral measure of a $G$-invariant Hilbert space ${\cal
H}_{[\mu]}$ of holomorphic functions on $\Xi$.

\Proof.  Since reproducing kernels are dominated by their values on the 
diagonal, we have to show only that

$$K_{[\mu]}(z,z)=\int_{\hat G_s} K_\lambda(z,z) \ d\mu(\lambda)$$ converges
uniformly on compact subsets of $\Xi$.  But this is immediate from Lemma
5.3(ii) and (5.5).\qed

\bobheadline{\S 6 Applications to harmonic analysis on $G/K$}

We shall give a holomorphic extension of the heat kernel on $G/K$ associated to the Riemannian structure from \S 4. 
This can be done quite naturally from the view point of our holomorphic extensions of spherical functions and it 
illustrates some complex analytic methods in the harmonic analysis on $G/K$. 
For use in a future paper we use this to formulate an analog of the Bargmann-Segal transform (cf.\ 
[St99] for a discussion on compact symmetric spaces).  

\subheadline{Holomorphic extension of the heat kernel}

Let $G/K$ be the Riemannian symmetric space associated to $G$ and $K$. In \S 4 
we took $X_1,\ldots, X_n$ an orthonormal basis of $\p$ and considered the 
differential operator

$$\Delta =\sum_{j=1}^n L_{X_j}^2$$ on $G_\C/K_\C$. It is well known that the 
Laplace-Beltrami operator on $G/K$ is given by the restriction of $\Delta$ to 
functions on $G/K$.  We remark
that with this convention $\Delta$ is a negative operator.  Denote also
by $\Delta$ the self-adjoint extension of $\Delta\res_{C_c^\infty(G/K)}$.
Fix a left $G$-invariant positive measure on $G/K$, say $\mu_{G/K}$ ,  and let
$L^2(G/K)$ be the corresponding $L^2$-space.  We will denote the identity coset 
in $G/K$ by $x_0$.

\par Consider the {\it heat equation} for $\Delta$:

$${\partial \over \partial t} u =\Delta u \qquad (u\in C^\infty(\R^+\times
G/K)).\leqno(6.1)$$ We write $k_t(x)$, $t\in \R^+$, $x\in G/K$, for the
heat kernel normalized by $$\lim_{t\downarrow 0} \int_{G/K} f(x)\
k_t(x) \ d\mu_{G/K}(x)=f(x_0) \leqno(6.2)$$ for $f\in C_c(G/K)$.  Recall
from [G68] the formula

$$k_t(x)={c\over |{\cal W}|} \int_{i\a^*} e^{-t(|\lambda|^2 +|\rho|^2)} \
\phi_\lambda(x) \ {d\lambda\over |c(\lambda)|^2}\qquad (x\in G/K)\
.\leqno(6.3)$$ Here $|\cdot|$ denotes the inner product on $\a_\C^*$ obtained
from the Cartan-Killing form, $c(\lambda)$ is the Harish-Chandra $c$-function
on $G/K$ while $c$ is a positive constant.  Substituting the integral formula 
for spherical functions into (6.3) we get

$$k_t(x)={c\over |{\cal W}|} \int_{i\a^*} \int_K e^{-t(|\lambda|^2
+|\rho|^2)} a(kx)^{\rho-\lambda} \ dk\ {d\lambda\over
|c(\lambda)|^2}\qquad(x\in G/K) .\leqno(6.4)$$

\par For $f\in L^p(G/K)$, $1\leq p< \infty$, 

$$u_f(x,t)\:=(k_t*f)(x)=\int_{G/K} k_t(g^{-1}x)\ f(gK)\
d\mu_{G/K}(gK)$$ is a solution of (6.1). From the normalization condition
(6.2), $u_f$ satisfies the initial condition $u_f(0,x)=f(x)$. Rephrased in terms 
of the generator for the heat semigroup, it becomes $e^{t\Delta}f=k_t*f$ for all 
$f\in L^p(G/K)$,
$1\leq p<\infty$.

\Theorem 6.1.  Let $G/K$ be a Riemannian symmetric space and $k_t(x)$ its
heat kernel.  Then $k_t$ has an extension to a function on
$\R^+\times \Xi$ which is holomorphic on $\Xi$.

\Proof.  First recall from (1.1) that the Iwasasa projection $a\:  G/K\to A$
extends holomorphically to $a\:  \Xi \to T_\Omega$.  Hence the theorem will
follow from (6.4) provided we can show that for any compact subset $Q\subeq
\Xi$ there exists a constant $C_Q>0$ such that

$$(\forall \lambda\in i\a^*)\qquad \sup_{z\in Q} |a(z)^{\rho-\lambda}| \leq
C_Q e^{|\lambda|_\Omega}.\leqno(6.5)$$

\par Since $Q$ is compact, we have that $$\sup_{z\in Q}|\Re \log
a(z)|<\infty .$$ But then (6.5) is immediate from (1.1) and the theorem
follows.\qed

If $f$ is a function on $G/K$ which extends holomorphically to $\Xi$, we
write $f^\sim$ for the holomorphic extension.

In view of Theorem 6.1, the map $$H_t\:  C_c(G/K)\to {\cal O}(\Xi), \
f\mapsto H_t(f); \ H_t(f)(z)=\int_{G/K} k_t^\sim(g^{-1}z)\ f(gK)\
d\mu_{G/K}(gK)\leqno (6.6)$$ is a well-defined $G$-equivariant map.  
The map $H_t$ is a  G-equivariant version of the familar
Bargmann-Segal  transform which we refer to as the  {\it heat kernel transform}. 
Notice that we have

$$H_t(f)=(k_t* f)^\sim$$ for all $f\in C_c(G/K)$.

\subheadline{The Fourier transform on $G/K$ and the Paley-Wiener space}

Helgason [He94, Ch.\ III, \S 1] defines the {\it Fourier-transform} on $G/K$ by

$${\cal F}\:  L^2(G/K)\to L^2(K/M\times i\a^*, \ d\mu_{K/M}\otimes
{d\lambda\over |c(\lambda)|^2})$$ $$({\cal F}f)(kM,\lambda)=\int_{G/K} f(gK)\
a(k^{-1}g)^{\rho-\lambda}\ d\mu_{G/K}(gK).$$ This is an isomorphism of Hilbert
spaces with inverse $${\cal F}^{-1}(F)(gK)=\int_{K/M\times i\a^*} F(kM,
\lambda)\ a(k^{-1}g)^{\rho+\lambda} \ d\mu_{K/M}(kM)\ {d\lambda\over
|c(\lambda)|^2}.$$ Define the {\it Paley-Wiener space} by

$${\rm PW}(G/K)={\cal F}^{-1}(C^\omega(K/M)\otimes C_c^\infty(i\a^*)).$$ Then 
this is a dense subspace of $L^2(G/K)$.

\Lemma 6.2.  All functions in ${\rm PW}(G/K)$ extend holomorphically to
$\Xi$, i.e., we have a $G$-equivariant injective mapping

$${\rm PW}(G/K)\to {\cal O}(\Xi), \ \ f\mapsto f^\sim.$$

\Proof.  Given the definition of the Paley-Wiener space and the Fourier
transform as explained above, the proof reduces to the estimate (6.5) which
was established in the proof of Theorem 6.1.\qed

By the definition of the heat kernel we have

$${\cal F}(k_t)(\lambda, kM)=e^{-t(|\lambda|^2 +|\rho|^2)}.\leqno(6.7)$$

\Lemma 6.3.  Let $t>0$.  Then the following assertions hold:  \item{(i)}
$k_t * {\rm PW}(G/K)\subeq {\rm PW}(G/K)$;  \item{(ii)} If $f\in {\rm
PW}(G/K)$, then $k_t*f$ extends holomorphically to $\Xi$.  In particular,
the heat kernel transform is well defined on ${\rm PW}(G/K)$ and we obtain
a $G$-equivariant map

$$H_t\:  {\rm PW}(G/K)\to {\cal O}(\Xi), \ \ f\mapsto (k_t*f)^\sim.$$

\Proof.  In view of Lemma 6.2, (ii) follows from (i).  Thus it suffices
to prove (i).  For that let $f\in {\rm PW}(G/K)$.  By definition of ${\rm
PW}(G/K)$ we have $f={\cal F}^{-1}(F)$ for some $F\in C^\omega(K/M)\otimes
C_c^\infty(i\a^*)$.  From (6.7) we obtain for all $g\in G$ that

$$\eqalign{(k_t*f)(gK)&={\cal F}^{-1}\big({\cal F}(k_t*f)\big) ={\cal
F}^{-1}\big({\cal F}(k_t){\cal F}(f)\big)= {\cal F}^{-1}\big({\cal
F}(k_t)F)\cr &=\int_{K/M}\int_{i\a^*} e^{-t(|\lambda|^2 +|\rho|^2)}\ F(kM,
\lambda) \ a(k^{-1}g)^{\rho+\lambda} \ {d\lambda\over |c(\lambda)|^2} \
d\mu_{K/M}(kM). \cr}\leqno(6.8)$$ Then the assertion in (i) is immediate from 
this.
\qed

\par Consider ${\cal O}(\Xi)$ as a Fr\'echet space with the topology
of compact convergence.  From (6.5) and (6.8) one can show
 that the heat kernel transform $H_t$ is well-defined on $L^2(G/K)$ for 
$t>0$. Moreover,

$$H_t\:  L^2(G/K)\to {\cal O}(\Xi), \ \ f\mapsto (k_t*f)^\sim$$ is an injective 
$G$-equivariant continuous mapping. Consequently, ${\cal 
G}_t(\Xi)\:=H_t(L^2(G/K)$ equip\-ped with the topology of
$L^2(G/K)$ is a $G$-invariant Hilbert space of holomorphic functions on $\Xi$
(cf.\ Section 5), hence has a reproducing kernel $${\cal
K}^t\:  \Xi\times \Xi\to \C, \ \ (z,w)\mapsto {\cal K}^t(z,w)$$ holomorphic
in the first and antiholomorphic in the second variable.  Using the technique to 
prove Theorem 5.1 one readily obtains 
\Theorem 6.4. The reproducing kernel for ${\cal G}_t(\Xi), {\cal K}^t(z,w)$ is 
given by
$${\cal K}^t(z,w)=\int_{G/K} k_t(g^{-1} z)\ \oline {k_t(g^{-1}w)} \
d\mu_{G/K}(gK) \qquad (z,w\in \Xi).$$ Moreover, 

$${\cal K}^t(z,w)=\int_{i\a^*} K_\lambda(z,w) \ e^{-2t(|\lambda|^2
+|\rho|^2)}\ {d\lambda\over |c(\lambda)|^2}$$ uniformly in $z,w\in \Xi$.
\qed

\bobheadline{\S 7 Hermitian symmetric subdomain of $\Xi$}

In addition to the intrinsic importance of the domain $\Xi$, it also contains 
some natural subdomains whose properties we think should be investigated. An 
introduction to these subdomains will be the topic of these last two sections. 
In this section we shall show that $\Xi$ coming from a classical group $G$ 
contains a  large $G$-invariant subdomain, $\Xi_0\subeq  \Xi$, that is 
$G$-biholomorphic to a Hermitian symmetric space of tube type. In some cases 
$\Xi_0 = \Xi$, and these we identify. For such groups $G$ one has the intriguing 
situation that the maximal Grauert domain of a Riemannian symmetric space of 
non-compact type, $G/K$,  is bi-holomorphic though not isometric to a Hermitian 
symmetric space. Hence, from Proposition 1.3, the harmonic analysis on $G/K$ 
extends holomorphically to this Hermitian symmetric space. These results were 
explicitly mentioned in [KS04] and the point of this section is to present the 
detailed argument.

\par The setup for this structure originates, we believe, with the paper by 
Nagano ([N65]) and subsequently has had many incarnations as, e.g., symmetric 
R-spaces, real forms of Hermitian symmetric spaces, real forms of simple Jordan 
triple systems, causal symmetric spaces, etc. We shall move freely between these 
different equivalent versions as necessary, but shall begin from the viewpoint 
of Jordan theory.

\par We consider first Riemannian symmetric spaces $G/K$ associated to Euclidean 
Jordan algebras in order to illustrate how $\Xi_0$ arises, but more importantly 
because this will lead to another domain in \S 8 seemingly more related to 
algebraic geometry.

\subheadline{$G/K$ associated to Euclidean Jordan algebras}

\par   We recall some basic
facts about Jordan algebras, referring the reader to [FK94] for a particularly 
nice presentation of the details. Let $V$ be a Euclidean Jordan algebra.  For 
$x\in V$ define $L(x)\in
\End(V)$ by $L(x)y=xy$, $y\in V$.  We denote by $e$ the identity element of
$V$.  The symmetric cone $W\subeq V$ can be defined by
$$W=\Int\{ x^2\:  x\in V\},$$ where $\Int\{\cdot\}$ denotes the topological 
interior of
$\{\cdot\}$.  Let $G=\Aut_0(W)$ be the identity component of the
automorphism group of $W$.  Then $G$ is a
reductive subgroup of $\Gl(V)$.  The isotropy group at $e$ $$K=G_e=\{
g\in G\:  g(e)=e\}$$ is a maximal compact subgroup of $G$, and the map
$$G/K\to W, \ \ gK\mapsto g(e)$$ is a homeomorphism.  \par Complexify $V$ to 
$V_\C
=V\oplus iV$ and consider the tube domain $$T_W=V+i
W\subeq V_\C.$$  The identity component of the complex
automorphism group of $T_W$ we denote by $G^h$.  Set
$K^h=G^h_{ie}$, the stabilizer of $ie\in T_W$.  Then $K^h$ is a maximal
compact subgroup of $G^h$ and contains $K$. The map $$G^h/K^h\to T_W, \ \
gK^h\mapsto g(ie)$$ is a biholomorphism of the Hermitian symmetric space
$G^h/K^h$ onto $T_W$ giving an explicit realization of $G^h/K^h$ as a tube 
domain.  Clearly one has a natural inclusion  of the real symmetric space $G/K$ 
into the Hermitian space $G^h/K^h$.

\par For the convenience of the reader we list the irreducible Euclidean
Jordan algebras and their associated groups $G$ and $G^h$ (cf.\ [FK94, p.\
213]).

\msk \centerline{{\bf Table I}} \centerline{{\bf Irreducible Euclidean Jordan
algebras}}

$$\vbox{\tabskip=0pt\offinterlineskip \def\tablerule{\noalign{\hrule}}
\halign{\strut#&\vrule#\tabskip=1em plus2em& \hfil#\hfil&\vrule#&\hfil#\hfil
&\vrule\hfil#\hfil& \hfil#\hfil&\vrule# \tabskip=0pt\cr\tablerule
&&\omit\hidewidth $V$\hidewidth&& \omit\hidewidth $G$\hidewidth&&
\omit\hidewidth $G^h$\hidewidth& \cr\tablerule && $\Symm(n,\R)$ &&
$\Sl(n,\R)\times \R^+$ && $\Sp(n,\R)$ &\cr\tablerule && $\Herm(n,\C)$ &&
$\Sl(n,\C)\times \R^+$ && $\SU(n,n)$ & \cr\tablerule && $\Herm(n,\H)$ &&
$\Sl(n,\H)\times \R^+$ && $\SO^*(4n)$ & \cr\tablerule && $\R\times\R^n$ &&
$\SO(1,n)\times\R^+$ && $\SO(2,n+1)$ &\cr\tablerule && $\Herm(3,{\Bbb O})$ &&
$E_{6(-26)}\times \R^+$ && $E_{7(-25)}$ & \cr\tablerule}} $$

\par Next let $c_1,\ldots, c_l$ be a Jordan frame (cf.\ [FK94, p.\ 44]) of $V$
and set $V^0=\bigoplus_{j=1}^l \R c_j$.  Then $K(V^0)=V$ (cf.\
[FK94, Cor.\ IV.2.7]).

\Lemma 7.1. $T_W=G\cdot ((V^0)_\C \cap T_W)$.

\Proof.  The inclusion $''\supeq''$ is clear.  Conversely, notice that
$$((V^0)_\C \cap T_W)=V^0+ i(W\cap V^0).$$ So let $z=x+iy\in T_W$, $x\in
V$, $y\in W$.  Then we find a $g\in G$ such that $g(y)=e$ and thus
$g(z)=g(x)+ie$.  Since $V=K(V^0)$, we find a $k\in K$ such that $k(g(x))\in
V^0$.  As $k(e)=e$ we get $(kg)(z)=(kg)(x)+ie\in ((V^0)_\C \cap T_W)$,
providing the other containment.  \qed

The choice $\a=\bigoplus_{j=1}^l \R L(c_j)$ defines a maximal abelian subspace 
orthogonal to $\k$.  As can be seen from Table I, the root system
$\Sigma=\Sigma(\g,\a)$ is classical and of type $A_{l-1}$.  If we define
$\eps_j\in \a^*$ by $\eps_j(L(c_i))=\delta_{ij}$, then we have $$\Sigma=\{
{1\over 2}(\eps_i-\eps_j)\:  i\neq j\}$$ (cf.\ [FK94, Prop.\ VI.3.3]).
Thus $$\Omega=\{ x=\sum_{j=1}^l x_j L(c_j)\:  x_j\in \R, \ |x_i-x_j|<\pi\}
.$$ In particular, $$\Omega_0\: =\bigoplus_{j=1}^l ]-{\pi\over 2},
{\pi\over 2}[ L(c_j)\subeq \Omega.$$  From Lemma 1.4 in [KS04] it follows that 
the domain
$$\Xi_0\: = G \exp(i\Omega_0)K_\C/ K_\C \subeq G_\C/ K_\C$$ is a $G$-invariant
open subdomain of $\Xi$. The interest in $\Xi_0$ is provided by the next result.

\Theorem 7.2.  The map $$\Xi_0\to T_W, \ \ gK_\C\mapsto g(ie)$$ is a
biholomorphism, i.e. $\Xi_0$ is biholomorphic to a Hermitian symmetric space.

\Proof.  In view of Lemma 7.1 it suffices to show that
$$T_\Omega(ie)=((V^0)_\C \cap T_W).$$ But since $((V^0)_\C \cap
T_W)=\bigoplus_{j=1}^l (\R+i\R^+)c_j$, this is immediate from the definition
of $T_\Omega$.\qed

\subheadline{Compactly causal symmetric spaces}

In fact, the characterization of $\Xi_0$ in the previous subsection will be seen 
to be a special case of results in this section. For the more general results 
one could use Jordan triple systems. Instead, we find it easier to use compactly 
causal symmetric spaces $S/G$.

\par Let $\s$ be a semisimple real Lie algebra
and $\theta$ a Cartan involution on $\s$ with Cartan decomposition
$\s=\uu\oplus \p_*$.  Let $\tau \: \s\to\s$ be an involution on $\s$ which, as
we may, commutes with $\theta$.  Write $\s=\g\oplus\q$ for the
$\tau$-eigenspace decomposition corresponding to the $\tau$-eigenvalues $+1$
and $-1$.  Set $\k=\g\cap\uu$ and $\p=\g \cap \p_*$. Then
$\g=\k\oplus \p$ is a Cartan decomposition of $\g$.

\par The symmetric Lie algebra $(\s,\tau)$ is called {\it irreducible} if the
only $\tau$-invariant ideals of $\s$ are $\{0\}$ and $\s$.  An irreducible
semisimple symmetric Lie algebra $(\s,\tau)$ is called {\it compactly causal} 
if, in the notation above, $\z(\uu)\cap\q\neq \{0\}$. \par For the rest of the 
section $\s$ will denote a compactly causal Lie algebra. For the convenience of 
the reader we recall the various types of compactly causal Lie algebras. 

\Remark 7.3.  (a) Suppose that $\s$ is simple. Then $\z(\uu)\neq \{0\}$ implies 
$\s$ is Hermitian, i.e. the symmetric space associated to $\s$ is Hermitian.  In 
this case we rename $\s$ to be $\g^h$. The compactly causal symmetric pairs 
$(\g^h,\g)$ are in Table II (cf.\ [H\'O96, Th.\ 3.2.8]).  \msk 
\centerline{{\bf Table II}}
\centerline{{\bf Compactly causal pairs $(\g^h,\g)$}} \centerline{{\bf for 
$\g^h$
simple}}

$$\vbox{\tabskip=0pt\offinterlineskip \def\tablerule{\noalign{\hrule}}
\halign{\strut#&\vrule#\tabskip=1em plus2em& \hfil#\hfil&\vrule#&\hfil#\hfil
&\vrule#\tabskip=0pt\cr\tablerule &&\omit\hidewidth $\g$\hidewidth&&
\omit\hidewidth $\g^h$\hidewidth & \cr\tablerule && $\sL(n,\R)\oplus \R$ &&
$\sp(n,\R)$ & \cr\tablerule && $\sL(n,\C)\oplus \R$ && $\su(n,n)$ &
\cr\tablerule && $\sL(n,\H)\oplus \R$ && $\so^*(4n)$ & \cr\tablerule &&
$\so(1,n)\oplus\R$ && $\so(2,n+1)$ & \cr\tablerule && $\e_{6(-26)}\oplus \R$
&& $\e_{7(-25)}$ & \cr\tablerule && $\so(p,q)$ && $\su(p,q)$ & \cr\tablerule
&& $\sp(p,q)$ && $\su(2p, 2q)$ & \cr\tablerule && $\so(n,\C)$ && $\so^*(2n)$
& \cr\tablerule && $\sp(n,\C)$ && $\sp(2n,\R)$ & \cr\tablerule &&
$\so(p,1)\oplus \so(q,1)$ &&$\so(2,p+q)$ & \cr\tablerule && $\sp(2,2)$ &&
$\e_{6(-14)}$ & \cr\tablerule && $\f_{4(-20)}$ && $\e_{6(-14)}$ &
\cr\tablerule && $\su^*(8)$ && $\e_{7(-25)}$ & \cr\tablerule}}$$

Notice that the first 5 listings in Table II are those $\g$ which arise as 
automorphism groups of the cone in a Euclidean Jordan
algebra. One calls such $(\g^h,\g)$ compactly causal symmetric pairs of {\it 
Cayley
type}.  \par\nin (b) Suppose that $\s$ is not simple. Then it turns out that 
$\s=\g\oplus\g$ with $\tau$
the flip involution $\tau(X,Y)=(Y,X)$ and $\g\cong\g^h$ a simple Hermitian Lie 
algebra (cf.\
[H\'O96, Lemma 1.3.7(2)]). We shall refer to this as the "group case". Of 
course the familiar simple Hermitian Lie algebras are, up to
isomorphism,

$$\su(p,q)\quad \so^*(2n)\quad \sp(n,\R)\quad \so(2,n)\quad \e_{6(-14)}\quad
\e_{7(-25)}.$$

\par\nin (c) From the tables in (a) and (b) one can observe the folklore fact 
that up to possibly a direct summand of $\R$, every classical simple real Lie 
algebra is the
fixed point algebra $\g$ in a compactly causal symmetric Lie algebra
$(\s,\tau)$, or said differently, every Riemannian symmetric space of noncompact 
type, $G/K$, with $G$ a classical group is the fixed point set of an 
anti-holomorphic involution on some Hermitian symmetric space. \qed

To generalize the previous subsection we need to define the analog of $\Omega_0$ 
and then $\Xi_0$. Unavoidably this will entail a little structure theory.
Let $\t\subeq \uu$ be a $\tau$-stable compact Cartan subalgebra of $\s$ and
$\Delta=\Delta(\g_\C,\t_\C)$ the associated root system.  Let
$X_0\in\z(\uu)\cap\q$ be normalized by $\Spec (\ad X_0)=\{ -i, 0,
i\}$. We follow the classic treatment of Hermitian symmetric spaces by 
Harish-Chandra. Choose a positive system $\Delta^+$ of $\Delta$ such that
$\Delta^+\subeq \{ \alpha\in \Delta\:  i\alpha(X_0)\geq 0\}$. The set of 
non-compact roots in $\Delta$, i.e., the set of
roots $\alpha$ for which $ \s_\C^\alpha\subeq \p_{*\C}$ is denoted $\Delta_n$, 
and, with$\p^{\pm}=\bigoplus_{\alpha\in \Delta_n^\pm} \s_\C^\alpha$, we have
$\p_{*\C}=\p^+\oplus\p^-$. Then $$\s_\C=\p^-\oplus\uu_\C\oplus\p^+.$$

\par Take $\Gamma=\{\gamma_1, \ldots, \gamma_r\}\subeq \Delta_n^+$ a
maximal set of long strongly orthogonal roots. By [H\'O96, Lemma
4.1.7] we may assume that $\Gamma$ is $-\tau$-invariant.  Let $\a_*\subeq
\p_*$ be the maximal abelian subspace of $\p_*$ which is constructed from
$\Gamma$.  Since $\Gamma$ is $-\tau$-invariant, $\a_*$ is
$\tau$-invariant. Hence $\a_*=\a\oplus\f$ with $\a=\a_*\cap\g$ and
$\f=\a_*\cap\q$.  Now the crucial observation is that $\a$ is maximal abelian in 
$\p$ by [H\'O96,
Lemma 4.1.9] so we define restricted root systems $\Sigma=\Sigma(\g,\a)$ and
$\hat\Sigma=\hat\Sigma( \s,\a)$.  Let $\sigma_1,\ldots, \sigma_s$ be the set
of restrictions of $\Gamma\circ C$ to $\a$ where $C\:  \s_\C\to\s_\C$ is the
Cayley transform with $C^{-1}(\a)\subeq i(\t\cap\g)$.  Now it turns out that
either $r=s$ (as for Cayley type spaces $(\g^h,\g)$) or $r=2s$
(as in the case $(\s,\g)=(\g^h\oplus\g^h, \g^h)$).

\par From various results of \'Olafsson one concludes that the
restricted root system $\hat\Sigma$ is either of type $C_s$ or $BC_s$, (for a 
proof cf.\ [K01, Th.\ 4.4]). Then

$$\hat \Sigma=\{ {1\over 2}(\pm\sigma_i\pm \sigma_j)\:  1\leq i,j\leq s\}\bs
\{0\}$$ or $$\hat\Sigma=\{{1\over 2}(\pm\sigma_i\pm \sigma_j)\:  1\leq
i,j\leq s\}\bs \{0\}\cup\{\pm{1\over 2} \sigma_i\:  1\leq i\leq s\}.$$
We are ready to define a domain $\Omega_0$ in $\a$. Set
$$\Omega_0=\Omega(\hat\Sigma)=\{X\in\a\:  (\forall \alpha\in \hat\Sigma)\
|\alpha(X)|<{\pi\over 2}\}.$$ Since$\Sigma\subeq \hat\Sigma$ it is clear that 
$\Omega_0\subeq \Omega = \Omega(\Sigma)$. Also, if we set $T_{\Omega_0}= 
A\exp(i\Omega_0)$ then $T_{\Omega_0}\subeq T_\Omega$.

\Lemma 7.4. $$\Omega_0=\{ X\in \a\:  (\forall \sigma_i)\
|\sigma_i(X)|<{\pi\over 2}\}.$$

\Proof.  Since $\hat\Sigma$ is of type $C_s$ or $BC_s$, this is immediate
from the definition of $\Omega_0$.\qed

\par To define $\Xi_0$ we remain with a compactly causal Lie algebra $\s$. As 
seen from examples (a) and (b) in Remark 7.3 $\s$ will be either a simple 
Hermitian algebra or the product of two such. Let $S$ be a connected Lie group 
with Lie algebra $\s$ and let $G,K,U$ the analytic subgroups of $S$ with Lie 
algebras $\g,\k$ and $\uu$.  We
will assume that $S$ is contained in its complexification $S_\C$.  Let
$P^\pm$ be the analytic subgroup of $S_\C$ corresponding to $\p^\pm$.

Write ${\cal D}\subeq \p^+$ for the Harish-Chandra realization of $S/U$
as a bounded symmetric domain.  It is convenient for us to consider ${\cal
D}$ as an open subset in the flag manifold $S_\C/ U_\C P^-$.  Hence $0\in
{\cal D}$ becomes identified with the identity coset in $S_\C/U_\C P^-$.

We define the domain $$\Xi_0= G\exp(i\Omega_0)K_\C/ K_\C.$$ Since 
$\Omega_0\subeq\Omega$, $\Xi_0$ is contained in $\Xi$ and by the result in 
[KS04] is a $G$-invariant open domain. Clearly, $\Xi_0$ contains 
$T_{\Omega_0}$.

\Theorem 7.5.  Let $S/G$ be a compactly causal symmetric space and let ${\cal
D}$ be the Harish-Chandra realization of the Hermitian symmetric space
$S/U$.  Then the map $$\Phi\:  \Xi_0\to {\cal D}, \ \ xK_\C\mapsto
x(0)$$ is a $G$-equivariant biholomorphism.  In particular, $\Xi_0$ is Stein.

\Proof.  First we show that $\Phi(T_{\Omega_0})\subeq {\cal D}$. So define 
elements $H_j\in \a$,
$1\leq j\leq s$, by $\sigma_k(H_j)=2\delta_{jk}$.  Recall
from Lemma 7.4 that $$\Omega_0=\bigoplus_{j=1}^s ]-{\pi\over 4}, {\pi\over
4}[ H_j.$$ In [H\'O96, Lemma 5.1.5] it is shown that there exists
pairwise orthogonal elements $E_1, \ldots, E_s$ of $\p^+\cap(\g \oplus i\q)$ (if
$r=s$, then $E_j\in{\s_\C}^{\gamma_j}$ while if $r=2s$, then, after a
possible renumbering of the elements in $\Gamma$, one has $E_j\in
{\s_\C}^{\gamma_{2j-1}}\oplus {\s_\C}^{\gamma_{2j}}$) such that for $x_j\in
\R$
$$\exp(\sum_{j=1}^s x_j H_j)(0)=\sum_{j=1}^s (\tanh x_j) E_j.$$  Since $\a 
+i\Omega_0=\sum_{j=1}^s \R H_j+i]-{\pi\over 4}, {\pi\over
4}[ H_j$ and the map $$\tanh\:  \R + i]-{\pi\over 4},
{\pi\over 4}[\to\{ z\in \C:  |z|<1\}, \ \ z\mapsto \tanh z$$ is a
biholomorphism, 
$$T_{\Omega_0}(0)=\{\sum_{j=1}^s z_jE_j\:  z_j\in\C,\ |z_j|<1\}.$$ But then
$T_\Omega(0)\subeq {\cal D}$.

\par We still have to show that $\Phi$ is onto.  From [S84, Prop.\
7.1.2] the map $$\p \times (\p_*\cap\q)\to S/U, \ \ (X,Y)\mapsto
\exp(X)\exp(Y)U$$ is a diffeomorphism (here $\p=\p_*\cap \g$).  Let
$\e$ be a maximal abelian subspaces of $\p_*\cap \q$.  All such$\e$ are 
conjugate under $\Ad(K)$. Hence $${\cal D}=G\exp(\e)(0).$$ So it would be 
enough to show there exists an $\e$ with $\exp(\e)(0)\subeq T_\Omega(0)$. Recall
the element $X_0\in \z(\uu)\cap\q$ and define an automorphism $J$ of $\s$
by $J=e^{\ad {\pi\over 2}X_0}$. Of course
$J$ also induces multiplication by $i$ on ${\cal D}$. A simple computation shows 
that
$J(\p)=\p_*\cap\q$ (cf.\ [\'O91, Lemma 1.4 (2)]).  In particular,
$\e=J(\a)$ defines a maximal abelian subspace $\e$ in $\p_*\cap\q$.   So we get 
$$\exp(\e)(0)=i\exp(\a)(0).$$ As $i\exp(\a)(0)\subeq
T_\Omega(0)$, the proof of the theorem is complete.  \qed

An immediate consequence of Theorem 7.5 in conjunction with Proposition 1.3 is
the following result. T. Kobayashi informs us that independently he and Faraut had 
also obtained this result in uncirculated personal notes.

\Theorem 7.6.  Assume that a Riemannian symmetric space $G/K$ is a totally
real form of a Hermitian symmetric space $S/U$ via a $G$-equivariant
embedding $G/K\into S/U$.  Then all eigenfunctions on
$G/K$ for the algebra of $G$-invariant differential operators $\D(G/K)$
extend holomorphically to $S/U$.\qed

\nin {\bf Example - the group case.}  Here we have $S=G\times G$ and $G=G^h$.  
In particular
$\uu=\k\times\k$ and $U=K\times K$.  Further $\t=\t_G\times\t_G$ with
$\t_G$ a compact Cartan algebra of $\k$ in $\g$.  The element
$X_0\in\z(u)\cap \q$ is such that $X_0=(X_0', -X_0')$ with $X_0'\in \z(\k)$
such that $\Spec\ad_\g(X_0')=\{-i, 0, i\}$.  Let $$\g_\C
=\p_G\oplus\k_\C\oplus\p_G^-$$ be the triangular decomposition of $\g_\C$
with respect to $X_0'$.  Then we have $$\p^+=\p_G^+\times \p_G^-\qquad
\hbox{and}\qquad \p^-=\p_G^-\times\p_G^+ .$$ Write ${\cal D}_G$ for the
Harish-Chandra realization of $G/K$ in $\p_G^+$ and $\oline {\cal D}_G$ for
the Harish-Chandra realization of $G/K$ in $\p_G^-$.  Then $${\cal D}={\cal
D}_G\times \oline{{\cal D}_G}.$$ Write ${\cal D}_G^{\rm opp}$ for ${\cal
D}_G$ equipped with the opposite complex structure.  Write $X\mapsto \oline
X$ for the complex conjugation in $\g_\C$ with respect to the real form $\g$.
Note that $$\oline {{\cal D}_G}\to {\cal D}_G^{\rm opp}, \ \ X\mapsto \oline
X$$ is a $G$-equivariant biholomorphism.  Denote by $0$ the origin in ${\cal
D}$.

Obviously we have $\hat\Sigma=\Sigma$ in these cases.  Thus $\Xi=\Xi_0$ and
Theorem 7.5 reads as:

\Theorem 7.7.  Assume that $G/K$ is a Hermitian symmetric space.  Then the
map $$\Phi\:  \Xi\to {\cal D}_G\times {\cal D}_G^{\rm opp}, \ \
xK_\C\mapsto (x(0), x(0))$$ is a $G$-equivariant biholomorphism.  In
particular, $\Xi$ is Stein.  \qed

\subheadline{The cases where $\Xi=\Xi_0$}

We have already seen that $\Xi=\Xi_0$ in the group case. Here we give some 
simple criteria to determine when $\Xi=\Xi_0$ in general. An independent and 
more differential geometric approach to this result is in [BHH03].

\Theorem 7.8.  The following are equivalent:  \item{(1)}
$\Xi=\Xi_0$;  \item{(2)} $\Omega_0=\Omega$;  \item{(3)} $\Sigma$ is of type
$C_s$ or $BC_s$ for $s\geq 2$ or $\Sigma=\hat\Sigma$ if $\Sigma$ has rank
one;  \item{(4)} $\rank_\R \g={1\over 2}\rank_\R \s$.

\Proof.  (1)$\iff$ (2) is just the definition.  \par\nin (2) $\iff$ (3):  In
the following we write $\Omega=\Omega(\Sigma)$ in order to make the
dependence clear on the root system.  For a classical root system $\Sigma$
one easily verifies the following facts:  \ssk

\item{$\bullet$} $\Omega(C_n)=\Omega(BC_n)$ for all $n\geq 1$.
\item{$\bullet$} $\Omega(B_n)=\Omega(D_n)$ for all $n\geq 3$.
\item{$\bullet$} $\Omega(C_n)\subsetneq \Omega(A_n)$ for $n\geq 2$.
\item{$\bullet$} $\Omega(C_n)\subsetneq \Omega(B_n)$ for $n\geq 3$.

\ssk\nin Recall that the root system $\hat \Sigma$ is of type $C_n$ or $BC_n$. A 
quick look at Table II tells us that the root system $\Sigma$ is always 
classical
for all $\g$ in question.  From that the equivalence (2) $\iff$
(3) is easily verified if $\rank\Sigma\geq 2$.  The various $\rank\Sigma=1$
cases one checks separately with Table II.  \par\nin (3) $\iff$ (4) This is
easily checked with Table II.\qed

{}From Theorem 7.8 we arrive at the following table:  \msk \centerline{{\bf
Table III}} \centerline{{\bf $\Xi=\Xi_0$ and $\s$ simple}}

$$\vbox{\tabskip=0pt\offinterlineskip \def\tablerule{\noalign{\hrule}}
\halign{\strut#&\vrule#\tabskip=1em plus2em& \hfil#\hfil&\vrule#&\hfil#\hfil
&\vrule#\tabskip=0pt\cr\tablerule &&\omit\hidewidth $\g$\hidewidth&&
\omit\hidewidth $\g^h =\s$\hidewidth & \cr\tablerule && $\sp(p,q)$ && $\su(2p,
2q)$ & \cr\tablerule && $\sp(n,\C)$ && $\sp(2n,\R)$ & \cr\tablerule &&
$\so(1,p)$ && $\so(2,p)$ & \cr\tablerule}}$$

\Remark 7.9.  (a) Suppose that $G/K$ is such that there is a Hermitian symmetric 
space $S/U$ containing $G/K$ as a totally real submanifold.  Using the
proof of Theorem 7.5 and the structure of $\Sigma$, one can show that not only 
the subdomain $\Xi_0$, biholomorphic to $S/U$, embeds in the projective variety
$S_\C/ U_\C P^-$, but so does $\Xi$ embed in it.  \par\nin (b) In [BHH03] the 
spaces
with $\Xi\neq \Xi_0$ are called rigid because, as they show, $\Aut_0(\Xi_0)=G$ 
if $\Xi\neq \Xi_0$. With hindsight and the observation in (a), this can be 
deduced easily using that the generic $S$-orbits in $S_\C/ U_\C P^-$ are open. 
Motivated by (a) we pose what seems to us also an interesting question. Consider 
the contraction semigroupof $\Xi$, namely $$\Gamma=\{ s\in S_\C\:  s\Xi\subeq 
\Xi\}.$$ Then we
ask if $\Gamma$ reduces to $G$ when$\Xi\neq \Xi_0$? 
 \par\nin (c) In an Appendix at the end of the paper we illustrate the main 
results in this section with two explicit compactly causal spaces.
 \qed

\bobheadline{\S 8  $\Xi^{1\over 2}$ - the ``square root'' of $\Xi$}

In this section we will describe another interesting domain associated to 
$G/K$'s coming from Euclidean Jordan
algebras $V$ which we call the square root
of $\Xi$, and denote by $\Xi^{1\over 2}$.  The square root domain $\Xi^{1\over 
2}$ is strictly smaller
than $\Xi$ and $\Xi_0$ but arises naturally as the maximal domain of
definition for the natural polarization of the metric on the symmetric cone
$W\simeq G/K$.

\subheadline{Characterization of $\Xi^{1\over 2}$}
\ssk Throughout this section $G/K$ is associated to a Euclidean Jordan
algebra $V$ as described earlier.  We will identify $G_\C/K_\C$ as a subset
of $V_\C$ by means of the orbit map $$G_\C/ K_\C \to V_\C , \ \ gK_\C
\mapsto g(e).  $$
We denote by $(\cdot|\cdot)$ the
Hermitian extension of the Jordan inner product on $V$ to $V_\C$.  Let $ \ 
z\mapsto \oline z$ be the complex conjugation with respect
to the real form $V$.  The technical aspects of the proofs will use the Pierce
decomposition of a Euclidean Jordan algebra (cf.\ [FK94, Th.\ IV.2.1]).
We recall it briefly. For every $\lambda\in \R$ and $x\in V$ we write 
$V(x,\lambda)$ for the$\lambda$-eigenspace of the symmetric operator $L(x)\in 
\End(V)$.  Set
$$V_{ij}=V(c_i, {1\over 2})\cap V(c_j, {1\over 2}) \qquad (1\leq i<j\leq n)$$
and $V_{ii}=\R c_i$ for $1\leq i\leq n$.  Then $$V=\bigoplus_{i\leq
j} V_{ij}\leqno(8.1)$$ and $$\eqalign{V_{ij}\cdot V_{ij}&\subeq V_{ii}
+V_{jj}\cr V_{ij}\cdot V_{jk}&\subeq V_{ik} \quad\hbox{if }\ i\neq k\cr
V_{ij}\cdot V_{kl}&=\{0\} \quad\hbox{if}\ \{i,j\}\cap \{
k,l\}=\eset.\cr}\leqno(8.2)$$ 
Here $V_{ij}\bot V_{kl}$ if $(i,j)\neq (k,l)$.

\par The {\it quadratic representation} (cf.\ [FK94, Ch.\ II, \S 3])
of the Jordan algebra $V_\C$ is defined by $$P(z)=2 L(z)^2 -
L(z^2)\in \End(V_\C) \qquad (z\in V_\C),$$ and its polarized version is defined 
by

$$P(z,w)=L(z)L(w)+ L(w)L(z) - L(zw) \qquad (z,w\in V_\C).$$ For this we have
the transformation property

$$P(gz, gw)=gP(z,w)g^t \qquad (\forall g\in G, \ z,w\in V_\C)\leqno (8.3)$$
(cf.\ [FK94, Prop.\ VIII.2.4]). 

\Lemma 8.1. For every $z\in V_\C$ the operator $P(z,\oline z)$ is Hermitian.
\Proof.  For $A\in\End (V_\C)$ let $\ A\mapsto A^*$ be the
conjugate transpose map.  Then $L(z)^*=L(\oline z)$ for all $z\in
V_\C$. From this and the definition of $P(z,w)$ the assertion follows.\qed

\par Set $\Omega^{1\over 2}={1\over 2}\Omega$ and define the {\it square
root} of $\Xi$ by $$\Xi^{1\over 2}= G \exp(i\Omega^{1\over 2}) K_\C /K_\C=
G \exp(i\Omega^{1\over 2})(e).$$ Again from [KS04] one has that $\Xi^{1\over 
2}$ is a $G$-invariant open subdomain in $G_\C/ K_\C$. From the definition of 
$\Xi_0$ one can see that $\Xi^{1\over 2}\subeq \Xi_0$. The next result gives an 
algebro-geometric characterization of $\Xi^{1\over 2}$.

\Proposition 8.2.
$$\Xi^{1\over 2}=\{ z\in GA_\C K_\C/ K_\C\subeq V_\C\:  P(z,\oline z)\ \hbox
{is positive definite}\}_0$$ where $\{\cdot\}_0$ denotes the connected
component of $\{\cdot\}$ which contains $e$.
\Proof. In view of the transformation property (8.3), we have to
show only that $$A\exp(i\Omega^{1\over 2})K_\C/ K_\C= \{ z\in A_\C(e)\subeq
(V^0)_\C\:  P(z,\oline z)\ \hbox {is positive definite}\}_0.$$ Fix
$z=\sum_{j=1}^l z_j c_j\in (V^0)_\C$, $z_j\in \C$. Let
$V=\bigoplus_{i\leq j}V_{ij}$ be the Pierce decomposition of $V$ (cf.\ (8.1),
(8.2)) and $V_\C=\bigoplus_{i\leq j}(V_{ij})_\C$ its complexification.  Let
$v\in (V_{ij})_\C$.  If $i\neq j$, then we have $L(c_k)v={1\over
2}(\delta_{ik}+\delta_{jk})v$, while for $i=j$ we have $L(c_k)v=\delta_{ik}
v$.  So the definition of $P(z,w)$ gives $$P(z,\oline
z)v=\cases{{1\over 2} (z_i\oline z_j + \oline z_i z_j)v & for $i\neq j$\cr
|z_i|^2v & for $i=j$.\cr}$$ 
Then $P(z,\oline z)$ is invertible if
and only if $$z_i\oline z_j + \oline z_i z_j\neq 0 \qquad (\forall 1\leq
i,j\leq l).\leqno(8.4)$$ \par Since $P(e)=\id$ is positive definite, it
follows from (i) and the continuity of the spectrum that $$\{ z\in A_\C(e)\:
P(z,\oline z)\ \hbox {is positive definite}\}_0 =\{ z\in A_\C(e)\:
P(z,\oline z)\ \hbox {is invertible}\}_0.$$ \par Suppose now that $z\in
A_\C K_\C / K_\C$.  Then $z=\exp(\sum_{j=1}^l w_j L(c_j))e$ for $w_j\in\C$
and we have $z_j=e^{w_j}$.  Hence (8.4) implies that $P(z,\oline z)$ is
invertible if and only if $$|\Im (w_i-w_j)|\not\in {\pi\over 2} +\Z\pi \qquad
(i\neq j).$$ In view of the definition of $\Omega^{1\over 2}$, this
completes the proof.\qed

\subheadline{The natural exhaustion function}

Proposition 8.2 implies in particular that the function $$\phi\:  \Xi^{1\over
2}\to \R, \ \ z\mapsto \log \det P(z,\oline z)^{-1}$$ is well defined and 
analytic.  Our goal is to show that $\phi$ is
plurisubharmonic.  \par Let $Z\in V_\C$ and $f$ a differentiable function
defined on some open subset in $V_\C$.  Then we set $$(\delta_Z f)(z)={d\over
dt}\Big|_{t=0}f(z+tZ).$$ For $Z\in V_\C$ we define the usual Cauchy-Riemann
operators  $$\partial_Z={1\over 2}(\delta_Z - i\delta_{iZ})\quad
\hbox{and}\quad \oline \partial_Z={1\over 2}(\delta_Z +i\delta_{iZ}).$$

\Lemma 8.3.  For all $Z_1, Z_2\in V_\C$ and $z\in \Xi^{1\over 2}$ we have
$$(\partial_{Z_1}\oline \partial_{Z_2}\phi)(z)= \tr\big[ P(z,\oline z)^{-1}
P(Z_1, \oline z) P(z,\oline z)^{-1} P(z, \oline Z_2)\big]- \tr
\big[P(z,\oline z)^{-1} P(Z_1, \oline Z_2)\big].$$

\Proof.  First notice that the definition of $\phi$ implies that

$$(\partial_{Z_1}\oline \partial_{Z_2}\phi)(z)=-{d\over dt}\Big|_{t=0}{d\over
ds}\Big|_{s=0} \log \det P(z+tZ_1,\oline z+ s\oline Z_2).$$

For all $t\in \R$ we have $$\eqalign{-{d\over ds}\Big|_{s=0}\log \det &
P(z+tZ_1,\oline z+ s\oline Z_2)=-{1\over \det P(z+tZ_1,\oline z)}\cdot\cr
&\cdot \det P(z+tZ_1,\oline z)\cdot \tr \big[P(z+tZ_1,\oline z)^{-1}{d\over
ds}\Big|_{s=0} P(z+tZ_1, \oline z+ s\oline Z_2)\big].\cr}$$ 
Since $${d\over ds}\Big|_{s=0} P(z+tZ_1, \oline z+ s\oline Z_2)= L(z+tZ_1)L(\oline
Z_2)+L(\oline Z_2)L(z+tZ_1)-L((z+tZ_1)\oline Z_2), $$ we obtain
$$\eqalign{-{d\over ds}\Big|_{s=0}\log \det & P(z+tZ_1,\oline z+ s\oline
Z_2)=- \tr\big[ P(z+tZ_1,\oline z)^{-1}( L(z+tZ_1)L(\oline Z_2)+\cr & +
L(\oline Z_2)L(z+tZ_1)-L((z+tZ_1)\oline Z_2))\big].\cr}$$ With that, finally we
get $$\eqalign{(\partial_{Z_1}\oline \partial_{Z_2}\phi)(z)
=&-{d\over dt}\Big|_{t=0} \tr\big[ P(z+tZ_1,\oline z)^{-1}( L(z+tZ_1)L(\oline
Z_2)+\cr &+L(\oline Z_2)L(z+tZ_1)-L((z+tZ_1)\oline Z_2))\big]\cr =&-\tr
\big[P(z,\oline z)^{-1}( L(Z_1)L(\oline Z_2)+ L(\oline Z_2)L(Z_1)-L(Z_1\oline
Z_2))+\cr &+\tr P(z,\oline z)^{-1} P(Z_1,\oline z) P(z,\oline z)^{-1}
P(z,\oline Z_2)\big],\cr}$$ concluding the proof of the lemma.\qed

\Lemma 8.4.  For all $g\in G$, $z\in \Xi^{1\over 2}$ and $Z\in V_\C$ we have
$$(\partial_{gZ}\oline \partial_{gZ}\phi)(gz)=(\partial_Z\oline
\partial_Z \phi)(z).$$

\Proof.  For all $g\in G$ we have by the transformation property (8.3) that

$$\eqalign{\phi(gz)& = -\log \det P(gz, g\oline z)=-\log \det gP(z, \oline
z)g^t\cr & =-\log \det P(z, \oline z)-\log \det gg^t=\phi(z)- \log \det
gg^t.\cr}$$ Therefore $$(\partial_{gZ}\oline
\partial_{gZ}\phi)(gz)=(\partial_Z\oline \partial_Z \phi\circ
g)(z)=\phi(z).\qeddis

\Lemma 8.5.  Let $V=\bigoplus_{1\leq i\leq j\leq l} V_{ij}$ be the Pierce
decomposition of a Euclidean Jordan algebra $V$.  Let $1\leq i,j,k,l,r,s\leq
l$ be integers such that $(i,j)\neq (k,l)$ and $i\neq j$ or $k\neq l$.  Then
the following assertions hold:  \item{(i)} $V_{ij}(V_{kl}\cdot V_{rs})\bot
V_{rs}$.  \item{(ii)} $(V_{ij}\cdot V_{kl})V_{rs}\bot V_{rs}$.

\Proof.  (i) By the invariance of the inner product on $V$ we have
$$\eqalign{ &V_{ij}(V_{kl}\cdot V_{rs})\bot V_{rs}\iff \cr & V_{kl}\cdot
V_{rs}\bot V_{ij}\cdot V_{rs}\iff \cr & V_{kl}\bot V_{rs}( V_{rs}\cdot
V_{ij}).\cr}$$ If $i\neq j$, then (i) follows from $V_{rs}( V_{rs}\cdot
V_{ij})\subeq V_{ij}$ which is seen as follows.  If $r=s$ or $(r,s)=(i,j)$,
then this is clear from (8.2).  So assume $r\neq s$ and $(r,s)\neq (i,j)$.
Also we may assume that $V_{rs}\cdot V_{ij}\neq \{0\}$.  Then $V_{rs}(
V_{rs}\cdot V_{ij})\not\perp V_{ij}$ again by the invariance of
$(\cdot|\cdot)$.  Since $(r,s)\neq (i,j)$ we obtain from (8.2) that
$V_{rs}\cdot V_{ij}=V_{uv}$ with $(u,v)\neq (r,s)$.  Hence $V_{rs}\cdot
V_{uv}\subeq V_{xy}$ and so $(x,y)=(i,j)$.  \par Let now $i=j$.  Then $k\neq
r$ by assumption.  It is enough to do the case $r=i$.  Then $$V_{rs}(
V_{rs}\cdot V_{ij})=V_{is}( V_{is}\cdot V_{ii})\subeq V_{is}\cdot V_{is}
\subeq V_{ii}+V_{ss}.$$ Now $V_{ii}+V_{ss}\bot V_{kl}$ since $k\neq l$
concluding the proof of (i).

\par\nin (ii) Again by the invariance of the inner product on $V$ we get
$$\eqalign{ &(V_{ij}\cdot V_{kl}) V_{rs}\bot V_{rs}\iff \cr & V_{ij}\cdot
V_{kl}\bot V_{rs}\cdot V_{rs}.\cr} $$ But $V_{rs}\cdot V_{rs}\subeq V^0$ and
$V_{ij}\cdot V_{kl}\subeq (V^0)^\bot$ by (8.2) and the assumptions in the Lemma. 
 \qed

The main result of this section is

\Theorem 8.6.  The function $$\phi\:  \Xi^{1\over 2}\to\R, \ \ z\mapsto \log
\det P(z,\oline z)^{-1}$$ is plurisubharmonic.

\Proof.  We have to show that $(\partial_Z\oline \partial_Z \phi)(z)\geq 0$
for all $Z\in V_\C$ and $z\in \Xi^{1\over 2}$.  In view of Lemma 8.4, we may
assume that $z\in \exp(i\Omega^{1\over2})K_\C/ K_\C \subeq (V^0)_\C$.  \par
For $Z_1, Z_2\in V_\C$ we define linear operators on $V_\C$ by $$A(Z_1,
\oline Z_2)=P(z,\oline z)^{-1} P(Z_1, \oline z)P(z,\oline z)^{-1} P(z,
\oline Z_2)$$ and $$B(Z_1,\oline Z_2)= P(z,\oline z)^{-1} P(Z_1, \oline
Z_2).$$ Then Lemma 8.3 implies that $$(\partial_{Z_1}\oline \partial_{Z_2}
\phi)(z)=\tr A(Z_1, \oline Z_2) -\tr B(Z_1, \oline Z_2).$$ \par Note that
$P(z,\oline z)^{-1}$ preserves all $(V_{ij})_\C$. In the proof of
Proposition 8.2 we have already shown that $$P(z,\oline
z)^{-1}v=\cases{{2\over (z_i\oline z_j + \oline z_i z_j)}v & for $v\in
V_{ij}$, $i\neq j$\cr |z_i|^{-2}v & for $v\in V_{ii}$.\cr}\leqno(8.5)$$

\ssk Then $V=V^0\oplus \bigoplus_{i\neq j} V_{ij}$ with
$V^0=\bigoplus_{j=1}^l \R c_j$.  Let now $Z=U+\sum_{i\neq j}Z_{ij}$ with
$U\in (V^0)_\C$ and $Z_{ij}\in (V_{ij})_\C$.  So

$$\eqalign{ (\partial_Z\oline \partial_Z \phi)(z)&=\tr (A(U,\oline
U)-B(U,\oline U)) + \sum_{i\neq j}\tr (A(Z_{ij},\oline
Z_{ij})-B(Z_{ij},\oline Z_{ij}))\cr & +\sum_{i\neq j} \tr (A(U,\oline Z_{ij})
+A(Z_{ij},\oline U) -2\Re B(U,Z_{ij})) -\sum_{ij\neq kl} \tr B(Z_{ij},
Z_{kl}) .\cr}$$

First we claim that $\tr A(U,\oline Z_{ij})=\tr A(Z_{ij},\oline U)=\tr
B(U,\oline Z_{ij})= \tr B(Z_{ij}, \oline Z_{kl})=0$.  In fact, this is an
immediate consequence of the complexified version of Lemma 8.5.  Hence we get
that $$(\partial_Z\partial_{\oline Z} \phi)(z)=\tr (A(U,\oline U)-B(U,\oline
U)) + \sum_{i\neq j}\tr (A(Z_{ij},\oline Z_{ij})-B(Z_{ij},\oline Z_{ij}))
.$$

We show separately that $\tr A(U,\oline U)-B(U,\oline U)\geq 0$ and $\tr
A(Z_{ij},\oline Z_{ij})-B(Z_{ij},\oline Z_{ij})\geq 0$.  \par We begin by
showing that $\tr A(U,\oline U)-B(U,\oline U)\geq 0$.  Write $U=\sum_{j=1}^l
u_j c_j$ for $u_j\in \C$.  Now $B(U,\oline U)=P(z,\oline z)^{-1} P(U,
\oline U)$ and

$$P(U,\oline U)v=\cases{|u_i|^2 v& if $v\in V_{ii}$,\cr {1\over 2} (u_i\oline
u_j+ u_j \oline u_i)v & if $v\in V_{ij}$, $i\neq j$.\cr}$$ From (8.5) we
obtain $$\tr B(U,\oline U)=\sum_{j=1}^n {|u_j|^2\over |z_j|^2}
+\sum_{i<j}{u_i \oline u_j + u_j\oline u_i\over z_i\oline z_j+ z_j \oline
z_i} \dim V_{ij}.$$ Next since $$P(U,\oline z)=L(U) L(\oline z) +L(\oline z)
L(U) - L(U\oline z)$$ we have $$P(U,\oline z)v=\cases{(u_i\oline
z_i)v& if $v\in V_{ii}$,\cr {1\over 2} (u_i\oline z_j + u_j\oline z_i)v & if
$v\in V_{ij}$, $i\neq j$\cr}.$$ Similarily one has $$P(z, \oline
U)v=\cases{(\oline u_i z_i)v& if $v\in (V_{ii})_\C$,\cr {1\over 2} (\oline
u_i z_j + \oline u_j z_i)v & if $v\in (V_{ij})_\C $, $i\neq j$.\cr}$$

Combined with (8.5) we get

$$\tr A(U,\oline U)=\sum_{j=1}^l {|u_j|^2\over|z_j|^2} + \sum_{i<j}
{|u_i\oline z_j + u_j \oline z_i|^2\over (z_i \oline z_j + z_j \oline z_i)^2}
\dim V_{ij}.$$ So $$\eqalign{\tr A(U,\oline U) -\tr B(U,\oline U)
&=\sum_{i<j} {\dim V_{ij}\over (z_i \oline z_j +z_j\oline z_i)^2}
\big(|u_i\oline z_j + u_j \oline z_i|^2 - (u_i\oline u_j+ u_j \oline u_i)
(z_i \oline z_j +z_j\oline z_i)\big)\cr &= \sum_{i<j}{\dim V_{ij}\over (z_i
\oline z_j +z_j\oline z_i)^2} |u_iz_j -u_j z_i|^2\geq 0\cr}\leqno(8.6)$$ as
we claimed.

\par Next we show that $\tr A(Z_{ij}, \oline Z_{ij})-B(Z_{ij},\oline
Z_{ij})\geq 0$ for $Z_{ij}\in (V_{ij})_\C$.  It is enough to do the case
$(i,j)=(1,2)$.  Set $e_{12}=c_1 +c_2$.  If $Z_1, Z_2\in (V_{12})_\C$, then
one has $$Z_1\cdot Z_2 ={1\over 2} (Z_1|\oline Z_2) e_{12}.$$ Fix $Z\in
(V_{12})_\C$.  We can normalize $Z$ such that $Z\oline Z= e_{12}$.  Fix $0\neq
Z\in (V_{12})_\C $ and set $$Z^{\bot_{12}}=\{ v\in (V_{12})_\C 
(v|Z)=0\}.$$ Similarily define $\oline Z^{\bot_{12}}$ We claim that
$ZQ=\oline Z Q =0$ for all $Q\in Z^{\bot_{12}}\cap \oline Z^{\bot_{12}}$.
In fact $ZQ, \oline Z Q\in (V^0)_\C $ by (8.2).  Let now $x\in (V^0)_\C $ be
arbitrary.  Then $x\oline Z=\lambda \oline Z$ for some $\lambda\in \C$ and so
$$(ZQ|x)=(Q|x\oline Z)=\oline \lambda(Q|\oline Z)=0$$ for all $x\in V^0$.
Thus $QZ=0$.  Similarily one shows that $Q\oline Z=0$, completing the proof
of our claim.  \par Polarizing [FK94, Lemma IV.2.2] gives

$$(L(Z)L(\oline Z)+L(\oline Z)L(Z))v={1\over 4} \|\ Z\|^2 v\quad\hbox {for
all} \ v\in V_{ij}, \ i\in \{1,2\}, \ j\geq 3.$$ Now $\|Z\|^2 =2$,
and since $Z\oline Z=e_{12}$ we have $P(Z,\oline Z)=L(Z)L(\oline Z) + L(\oline 
Z)
L(Z)-L(e_{12})$. So $$P(Z,\oline Z)v=\cases{c_2& if $v=c_1$,\cr c_1& if
$v=c_2$,\cr 0 & if $v\in \bigoplus_{j=3}^l \R c_j$\cr -v & if $v\in
Z^{\bot_{12}}\cap \oline Z^{\bot_{12}} $, \cr {1\over 2} (Z|\oline Z) \oline
Z & if $v=Z$,\cr {1\over 2} (\oline Z|Z) Z & if $v=\oline Z$,\cr 0 & if $v\in
V_{ij}$, $i\neq j$, $(i,j)\neq (1,2)$.\cr}$$ Using (8.5) we get
$$\tr B(Z,\oline Z)= {2\over z_1\oline z_2 + \oline z_1 z_2} ({1\over 2}
(Z|\oline Z)(\oline Z|Z) -\dim Z^{\bot_{12}}\cap \oline Z^{\bot_{12}}).$$
Again, $(Z|Z)=2$ together with the Cauchy-Schwarz inequality yields $$-\tr
B(Z,\oline Z)\geq {2\over z_1\oline z_2 + \oline z_1 z_2} (-2 +\dim
Z^{\bot_{12}}\cap \oline Z^{\bot_{12}}).$$ \par To compute $\tr A(Z,\oline
Z)$ first notice that $$P(z,\oline Z)=L(z)L(\oline Z)+L(\oline
Z)L(z)-{z_1+z_2\over 2}L(\oline Z).$$ A simple calculation now shows that
$$P(z,\oline Z)v=\cases{{z_1\over 2} \oline Z& if $v=c_1$,\cr {z_2\over
2}\oline Z & if $v=c_2$,\cr 0 & if $v\in V_{ij}$, $\{ i,j\}\cap \{
1,2\}=\eset$,\cr 0& if $v\in Z^{\bot_{12}}$, \cr z_1 c_1 + z_2 c_2 & if
$v=Z$,\cr z_j \oline Zv& if $v\in V_{ij}$, $i\in \{ 1,2\}$, $j\geq 3$.\cr}$$
Similarily we obtain that $$P(Z,\oline z)v=\cases{{\oline z_1\over 2} Z& if
$v=c_1$,\cr {\oline z_2\over 2}Z & if $v=c_2$,\cr 0 & if $v\in V_{ij}$, $\{
i,j\}\cap \{ 1,2\}=\eset$,\cr 0& if $v\in \oline Z^{\bot_{12}}$, \cr \oline
z_1 c_1 + \oline z_2 c_2&if $v=\oline Z$,\cr \oline z_j Zv& if $v\in V_{ij}$,
$i\in \{ 1,2\}$, $j\geq 3$.\cr}$$ With (8.5) this implies that $$A(Z,\oline
Z)v=\cases{{1\over z_1\oline z_2 + \oline z_1 z_2} c_1 + {1\over z_1\oline
z_2 + \oline z_1 z_2}{\oline z_1\over \oline z_2} c_2& if $v=c_1$,\cr {1\over
z_1\oline z_2 + \oline z_1 z_2} c_2 + {1\over z_1\oline z_2 + \oline z_1
z_2}{\oline z_2\over \oline z_1} c_1& if $v=c_2$,\cr 0& if $v\in V_{ij}$, $\{
i,j\}\cap \{ 1,2\}=\eset$,\cr 0& if $v\in Z^{\bot_{12}}$, \cr {2\over
z_1\oline z_2 + \oline z_1 z_2}v & if $v=Z$, \cr {|z_j|^2\over (z_1 \oline
z_j + z_j \oline z_1)(z_2 \oline z_j + z_j \oline z_2)}v & if $v\in V_{ij}$,
$\{ i,j\}\cap \{ 1,2\}=\eset$, $j\geq 3$.\cr}$$

Hence we get that $$\tr A(Z,\oline Z)={4\over z_1\oline z_2 + \oline z_1
z_2}+ \sum_{i\in \{ 1,2\}\atop j\geq3} {|z_j|^2\dim V_{ij}\over (z_1 \oline
z_j + z_j \oline z_1)(z_2 \oline z_j + z_j \oline z_2)} $$ and so $$\tr
A(Z,\oline Z) -B(Z,\oline Z)= {2\dim (Z^{\bot_{12}}\cap \oline
Z^{\bot_{12}})\over z_1\oline z_2 + \oline z_1 z_2} +\sum_{i\in \{ 1,2\}\atop
j\geq3} {|z_j|^2\dim V_{ij}\over (z_1 \oline z_j + z_j \oline z_1)(z_2 \oline
z_j + z_j \oline z_2)} >0.\leqno(8.7)$$ This concludes the proof of the
theorem.  \qed

\subheadline{ $\Xi^{1,{1\over2}}$ - the commutator square root domain}

The function $\phi\:  \Xi^{1\over 2}\to \R$ is almost strictly
plurisubharmonic, in the sense that it is strictly plurisubharmonic on a
codimension one subdomain.

\par Write $\g^1=[\g,\g]$ for the commutator subalgebra of $\g$. Then 
$\g=\g^1\times\R$.  With $\p^1=\g^1\cap \p$ we obtain a
Cartan decomposition of $\g^1$, $\g^1=\k\oplus\p^1$.  Denote by $G^1$ the
analytic subgroup of $G$ corresponding to $\g^1$ and by $G_\C^1$ its
complexification.  Then one has the canonical isomorphism $G_\C /K_\C \simeq
G_\C^1/ K_\C \times \C^*$ and the canonical embedding $$\iota\:  G_\C^1/ K_\C
\into G_\C/ K_\C , \ \ xK_\C\mapsto (xK_\C, 1).$$ Henceforth we identify
$G_\C^1/K_\C$ as a complex submanifold of $G_\C/ K_\C$ via the embedding
$\iota$.  \par Set $\Omega^1=\Omega\cap\g^1$ and notice that
$\Omega=\Omega^1\times\R$.  Then we define the {\it commutator square root
domain} by

$$\Xi^{1,{1\over2}}=G^1 \exp(i{1\over 2}\Omega^1) K_\C/ K_\C.$$ Clearly 
$$\Xi^{1\over 2}\simeq \Xi^{1, {1\over 2}}\times \C^*.$$ We shall consider 
$\Xi^{1,{1\over2}}$ as a submanifold of $\Xi^{1\over 2}$ via the
embedding $\iota$.  \par Set $V_\C^1=\{ z\in V_\C \:  \det z=1\}$.  Note that
$V_\C^1$ is a complex hypersurface in $V_\C$ and $\Xi^{1, {1\over2}}$ is open
in $V_\C^1$.

\par Consider the function

$$\phi^1\:  \Xi^{1,{1\over 2}}\to \R, \ \ z\mapsto \log \det P(z,\oline
z)^{-1},$$ the restriction of $\phi$ to
$\Xi^{1,{1\over 2}}$.

>From the proof of Theorem 8.6 we can obtain

\Theorem 8.7.  The function $$\phi^1\:  \Xi^{1,{1\over 2}}\to \R, \ \
z\mapsto \log \det P(z,\oline z)^{-1}$$ is strictly plurisubharmonic.

\Proof.  We have to show that $$(\partial_Z \oline \partial_Z
\phi)(z)>0\leqno(8.8)$$ for all $z\in \Xi^{1, {1\over 2}}$ and $Z\in T_z
\Xi^{1, {1\over 2}}$, $Z\neq 0$.  By Lemma 8.4 we may assume that $z\in
\exp(i{1\over 2}\Omega^1)(e)$, hence $z=\sum_{j=1}^l z_j c_j$ with $z_j\in
\C$.  Since $\Xi^{1, {1\over2}}$ is open in $V_\C^1$, we have
$$T_z\Xi^{1,{1\over 2}}=T_z V_\C^1.$$ Further $T_eV_\C^1=\{ u\in V_\C \:
\tr u=0\}$ and so

$$T_z\Xi^{1,{1\over 2}}= dL_z(e) \big(T_eV_\C^1=\{zu\:  u\in V_\C, \tr u=0\}
\big).$$ In particular, we obtain that

$$T_z\Xi^{1,{1\over 2}}=\bigoplus_{i\neq j} (V_{ij})_\C \oplus \big\{
\sum_{j=1}^l z_j u_j c_j\:u_j\in\C, \ \sum_{j=1}^l u_j=0\big\}.$$ Now (8.8)
follows from (8.6) (which is positive now if $U=\sum_{j=1}^lz_j u_j c_j\neq
0$ and $\sum_{j=1}^l u_j=0$) and (8.7).\qed

We give some applications of Theorem 8.6 to the complex analysis of
the square root domain $\Xi^{1\over 2}$.  In particular, we will show that
$\Xi^{1\over 2}$ is Stein and that there is a natural $G$-invariant K\"ahler
structure on $\Xi^{1,{1\over 2}}$.

\Theorem 8.8.  The square root domain $\Xi^{1\over 2}\subeq G_\C/K_\C$ is
Stein.

\Proof. In view of Grauert's solution of the Levi problem (cf.\ [H\"o73,
Th.\ 5.2.10]), it suffices to find a strictly plurisubharmonic exhaustion
function of $\Xi^{1\over 2}$.  Evidently $$\phi_0\:  \Xi^{1\over 2}\to \R^+,\
\ z\mapsto \|z\|^2$$ is strictly plurisubharmonic. Using the $\phi$ of Theorem 
8.6 we set
$\phi_1(z)\:=\exp \phi(z)$. Since
$\phi$ is plurisubharmonic, $\phi_1$ is plurisubharmonic.  Since $\phi_1$ is 
positive
$$\psi=\phi_0 +\phi_1$$ is a positive strictly
plurisubharmonic function on $\Xi^{1\over 2}$.  It remains to see that $\psi$
is proper.  For that let $(z_n)_{n\in\N}$ be a sequence which tends to
$\infty$ in $\Xi^{1\over 2}$.  Then we have either $z_n\to \infty$ in $V_\C$
or $z_n\to z_0\in \partial \Xi^{1\over 2}$.  In the first case $\phi_0$ blows
up, while in the second case $\phi_1$ blows up.  This concludes the proof of
the theorem.\qed

\subheadline{The K\"ahler structure}

\Proposition 8.9.  For $z\in \Xi^{1\over 2}$ the family of Hermitian forms
$H_z(\cdot,\cdot)$ on $T_z\Xi^{1\over 2}$ defined by

$$H_z(v,w)=(P(z,\oline z)^{-1}v|w) \qquad (v,w\in V_\C\simeq T_z
\Xi^{1\over 2})$$ defines a $G$-invariant Hermitian
metric on $\Xi^{1\over 2}$.  Moreover the associated Riemannian
structure $G_z(\cdot.\cdot)=\Re H_z(\cdot, \cdot)$ is
complete.

\Proof.  It follows from Proposition 8.2 and the transformation property
(8.3) that $H_z(\cdot,\cdot)$ defines a $G$-invariant Hermitian
structure on $\Xi^{1\over 2}$.  That the associate Riemannian metric is
complete follows from the fact that the metric blows up at the boundary.\qed

\Remark 8.10.  (a) If $z\in V_\C$, write $z=x+iy$ with $x,y\in V$.  Then
 a simple computation shows that $$P(z,\oline z)=P(x)+P(y).\leqno(8.9)$$ In
 particular, as a refinement of Proposition 8.2 we have the square root domain 
$\Xi^{1\over 2}$
 equivalently defined as

$$\Xi^{1\over 2}=\{ z=x+iy\in V_\C\:  P(x)+P(y)\in \End(V)\quad
\hbox{is positive definite} \}_0.$$

\par\nin (b) With $H_z(\cdot,\cdot)$ as in
Proposition 8.9 set $H_z=G_z + i\Omega_z$, $\Omega_z$ the imaginary part of 
$H_z$.  We ask whether the Hermitian metric $H_z(\cdot,\cdot)$ is K\"ahler. From
(8.9) we see that $P(z,\oline z)$ preserves the real form $V$.  Identify $V_\C$ 
with $V\times V$ via $V\times V\to V_\C, \
(x,y)\mapsto x+iy$.  Then in matrix notation we have
$$G_z=\pmatrix{P(z,\oline z)^{-1}\res_V & 0 \cr 0 & P(z,\oline
z)^{-1}\res_V\cr}\quad \hbox{and} \quad \Omega_z=\pmatrix{0 & P(z,\oline
z)^{-1}\res_V \cr -P(z,\oline z)^{-1} \res_V & 0\cr}.$$ If $V$ is
irreducible, then one can show that the Hermitian metric $H_z$ is K\"ahler if
and only if $V=\R$.

\par\nin (c) The $G$-invariant Riemannian metric on the cone $W$ is known 
[FK94, Th.\ III.5.3] to be given by $$g_x(v,w)=(P(x)^{-1}v|w)
\qquad (x\in W\ ; v,w\in V\simeq T_xW).$$ Since $P(z,z)=P(z)$, we see that the 
Hermitian metric $H_z$ is a polarization of the Riemannian
metric on $W$ to the square root domain $\Xi^{1\over 2}$. 

\par\nin (d) We recall from [FK94, p.  15-16 and Prop.\ III.4.3] that the
metric $(g_x)_{x\in W}$ is obtained from a potential function $$\psi\:  W\to
\R, \ \ x\mapsto {1\over 2}\log \det P(x)^{-1}$$ in such a way that
$$g_x(u,v)=(\partial_u\partial_v\psi)(x) \qquad (x\in W; \ u,v\in V\simeq T_x
W).$$ One might now expect that the Hermitian metric $H_z$ is obtained in a
similar manner through the potential function $$\phi\:  \Xi^{1\over 2} \to
\R, \ \ z\mapsto {1\over 2} \log \det P(z,\oline z)^{-1}$$ and the associated
Hermitian metric

$$h_z(u,v)=(\partial_u\oline \partial_v \phi)(z)\qquad (z\in\Xi^{1\over 2}; \
u,v\in V_\C\simeq T_z \Xi^{1\over 2}).$$ However, it turns out that
$h_z\neq H_z$.  The easiest example to observe this phenomenon is for
$V=\R$.  Then $W=\R^+$ and $G=\R^+$ is the group of scaling transformations.
Here we have $$P(x)=x^2,\quad \psi(x)=-\log x\quad \hbox{and}\quad
g_x(u,v)={1\over x^2} uv \qquad (x\in \R^+=W, u,v\in\R).$$ On the other hand
we have $\Xi^{1\over 2}=\C^*$ and $P(z,\oline z)=|z|^2$.  Hence
$\phi(z)=-\log |z|$.  Since $\phi$ is harmonic, we thus obtain $h_z=0$ and so
evidently $h_z\neq H_z$.\qed

Finally we construct the $G$-invariant K\"ahler structure on 
$\Xi^{1,{1\over2}}$.

\Theorem 8.11. Define $h_z \: T_z\Xi^{1,{1\over 2}}\times
T_z\Xi^{1,{1\over 2}}\to \C $ by $$ (Z_1, Z_2)\mapsto
(\partial_{Z_1}\oline\partial_{Z_2}\phi^1)(z)\qquad (z\in \Xi^{1, {1\over
2}}).$$ Then $h_z(\cdot,\cdot)$ is a $G^1$-invariant positive K\"ahler structure 
on 
$\Xi^{1,{1\over 2}}$ whose associated Riemannian structure
$(g_z)_{z\in\Xi^{1, {1\over 2}}} =(\Re h_z)_{z\in\Xi^{1, {1\over 2}}}$ is
complete.

\Proof.  This is immediate from Theorem 8.7 and Lemma 8.4.  The completeness
of the associated Riemannian structure follows from the fact that the metric
blows up at the boundary (cf.\ the formulas (8.6) and (8.7)).\qed

\Remark 8.12. (a) To us it seems an interesting question to determine if the 
K\"ahler metric$h_z$ is Einstein-K\"ahler. 
\par (b) Certainly there should be a generalization of the results in this 
section to the other groups $G$ whose Lie algebra appears in Table II.\qed

\bobheadline{Appendix - Examples of \S 7}
\bsk\nin {\bf The example of $G=\SO(p,q)$.}  We assume that $p\leq q$.
Choose $K=S(O(p)\times O(q))$ as a maximal compact subgroup of $G$.  For $G^h$ 
and $K^h$ we have $G^h=\SU(p,q)$ and $K^h=S(U(p)\times
U(q))$.  The involution $\tau$ on $G^h$ which has $G$ as fixed point group
is $$\tau\:  G^h\to G^h, \ g\mapsto \oline g$$ the complex conjugation,
i.e., if $g=x+iy$, $x,y \in M(n,\R)$, then $\oline g=x-iy$.

\par An appropriate maximal abelian subspace of $\p$ in $\g=\so(p,q)$ is
$$\a=\{ \pmatrix{ 0_{pp} & I_{t_1,\ldots, t_p}\cr 
I_{t_1,\ldots, t_p}^t & 0_{qq}\cr}\:  t_1,\ldots,t_p\in \R\},$$ where 
$$I_{t_1,\ldots,t_p}=\pmatrix{ & & & t_1\cr 0_{p,p-q}&  \addots & \cr 
& t_p & & \cr} \in M(p\times q;\R).$$ Let $n=p+q$.  For all $1\leq i,j\leq n$ 
define $E_{ij}\in M_n(\R)$ by$E_{ij}=(\delta_{ki}\delta_{jl})_{k,l}$.  Then 
$\a=\bigoplus_{j=1}^p \R e_j$ with $$e_j=E_{j, p+q+1-j}+ E_{p+q+1-j, j}.$$ 
Define $\eps_j\in \a_*$ by$\eps_j(e_k)=\delta_{jk}$. Then the root system 
$\Sigma =\Sigma(\a,\g)$ isgiven by $$\Sigma=\cases{ \{ \pm\eps_i\pm\eps_j\:  
1\leq i\neq j\leq p\}\cup\{ \pm\eps_i\:  1\leq i\leq p\} & for $1<p< q$,\cr 
\{\pm\eps_i\pm\eps_j\:  1\leq i\neq j\leq p\} & for $1<p=q$,\cr \{\pm \eps_1\} 
& for $p=1$.\cr}$$ In all cases one easily sees that

$$\Omega_0=\bigoplus_{j=1}^p]-{\pi\over 4}, {\pi\over 4}[ e_j$$ lies in
$\Omega$.

\par The Harish-Chandra realization of $G^h/K^h$ is given by $${\cal D}=\{
Z\in M(p\times q; \C)\:  I_q- Z^* Z>>0\}$$ where $G^h$ acts on ${\cal D}$ via
$$g(Z)=(AZ+B)(CZ+D)^{-1} \qquad (g=\pmatrix{ A & B\cr C & D\cr}\in G^h, Z\in
{\cal D}).$$ Define $${\cal D}_\R={\cal D}\cap M(p\times q; \R).$$ Note
that $\tau$ induces an invlution on ${\cal D}\simeq G^h/K^h$ which is given
by complex conjugation.  Hence $G/K$ defines a real subspace of $G^h/K^h$ and
since $G^h/G$ is compactly causal we even have $${\cal D}_\R\simeq G/K$$ by
[KN\'O97, Th.\ II.9].

Let $a=\exp(\sum_{j=1}^p z_j e_j)$ with $z_j\in \R + i]-{\pi\over 4},
{\pi\over 4}[$.  Then an easy calculation shows that

$$a(0)=\pmatrix{ & & & \tanh z_1\cr 0_{p,p-q}&  &\addots & \cr 
& \tanh z_p & & \cr}$$ and so

$$a(0)^* a(0)= \pmatrix{ |\tanh z_p|^2 & & \cr & \ddots & \cr
& & |\tanh z_1|^2 \cr}.$$

Hence $$\Phi(A\exp(i\Omega_0) K_\C/ K_\C)= \{ \pmatrix{ & & & z_1\cr
0_{p,p-q}&  &\addots & \cr & z_p & & \cr} 
\:  z_i\in \C, |z_i|<1\}\subeq {\cal D}.$$

The maximal abelian subspace $\e=J(a)$ of $\p_*\cap \q$ is given by

$$\e=\{ \pmatrix{ 0_{pp} & iI _{t_1,\ldots, t_p}\cr
-iI_{t_1,\ldots, t_p}^t& 0_{qq}\cr}\:  t_1,\ldots,t_n\in \R\}.$$ Then 
$$\exp(\e)(0)= \{ \pmatrix{& & & ix_1\cr 0_{p,p-q}&  &\addots & \cr 
& ix_p & & \cr} \:  x_j\in \R,|x_j|<1 \},$$ and so
$\exp(\e)(0)\subeq \Phi(A\exp(i\Omega_0) K_\C/ K_\C)$
which was the crucial step in the proof of Theorem 7.5.

\bsk\nin {\bf The example of $G=\SO(n,\C)$.}  We let $G^h=\SO^*(2n)$ realized
as

$$G^h=\{g\in \SU(n,n)  g^t \pmatrix { 0 & I\cr I& 0\cr} g =\pmatrix { 0 &
I\cr I& 0\cr}\}.$$ A maximal compact subgroup of $G^h$ is given
by $$K^h=\{ \pmatrix { A & 0\cr 0 & \oline A\cr}\:  A\in U(n)\} \simeq
U(n).$$ The Lie algebra $\g^h$ of $G^h$ is given by $$\g^h=\{ \pmatrix { A
& B\cr -\oline B & \oline A\cr}\:  A, B\in M(n,\C);\ A^*=-A,\ B^t=-B\}.$$
We equip $G^h$ with the involution $\tau\:  G^h\to G^h, \ g\mapsto \oline g$
of complex conjugation and set $G={G^h_0}^\tau$ for the connected component
of $\tau$-fixed points.  Then the Lie algebra $\g$ of $G$ is given by $$\g=\{
\pmatrix { A & B\cr -B & A\cr}\:  A, B\in M(n,\R);\ A^t=-A,\ B^t=-B\} .$$
Hence $\g\simeq \so(n,\C)$ and it is easy to check that $G\simeq
\SO(n,\C)$.  Then the Harish-Chandra realization of $G^h/K^h$ is given by
$${\cal D}=\{ Z\in M(n, \C)\:  Z^t=-Z,\ I- Z^* Z>>0\}$$ with $G^h$ acting
on ${\cal D}$ via $$g(Z)=(AZ+B)(CZ+D)^{-1} \qquad (g=\pmatrix{ A & B\cr C &
D\cr}\in G^h, Z\in {\cal D}).$$ Set ${\cal D}_\R={\cal D}\cap M(n,\R)$.
Now $\tau$ induces an involution on ${\cal D}\simeq G^h/K^h$ by complex 
conjugation.  Hence $G/K$ defines a real subspace of
$G^h/K^h$ and, since $G^h/G$ is compactly causal, we have ${\cal D}_\R\simeq
G/K$ by [KN\'O97, Th.\ II.9].

We have to distinguish the cases of $n$ even and odd.  Suppose first that
$n=2m$ is even.  Then as a maximal abelian subspace $\a$ of $\p\subeq \g$ we
choose

$$\a=\{ \pmatrix{ 0 & J_{t_1,\ldots, t_m}\cr -J_{t_1,\ldots, t_m} & 0\cr}\:
t_1,\ldots,t_m\in \R\},$$ where $$J_{t_1,\ldots, t_p}=\pmatrix{ \pmatrix {0
& t_1\cr -t_1 & 0\cr } & & \cr &\ddots& \cr & & \pmatrix {0 & t_m\cr -t_m &
0\cr}\cr}.$$ Then $\a=\bigoplus_{j=1}^m \R e_j$ with $$e_j=E_{2j-1, m+2j}-
E_{2j, m+2j-1}+ E_{m+2j, 2j-1} - E_{m+2j-1, 2j}$$

If $n=2m+1$, $m\geq 1$, is odd, then we choose

$$\a=\{ \pmatrix{ 0 & J'_{t_1,\ldots, t_m}\cr -J'_{t_1,\ldots, t_m} &
0\cr}\:  t_1,\ldots,t_m\in \R\}\ ,$$ where $$J'_{t_1,\ldots, t_p}=\pmatrix{
\pmatrix {0 & t_1\cr -t_1 & 0\cr } & & & \cr &\ddots& & \cr & & \pmatrix {0 &
t_m\cr -t_m & 0\cr} & \cr & & & 0\cr}.$$ In this case $\a=
\bigoplus_{j=1}^m \R e_j$.

Define $\eps_j\in \a^*$ by $\eps_j(e_k)=\delta_{jk}$.  Then the root system
$\Sigma =\Sigma(\a,\g)$ is given by

$$\Sigma=\cases{ \{ \pm\eps_i\pm\eps_j\:  1\leq i\neq j\leq m\} & for
$n=2m$,\cr \{ \pm\eps_i\pm\eps_j\:  1\leq i\neq j\leq m\} \cup\{ \pm\eps_i\:
1\leq i\leq m\}& for $n=2m+1$.\cr }$$

In both cases we see that the set

$$\Omega_0=\bigoplus_{j=1}^m]-{\pi\over 4}, {\pi\over 4}[ e_j.$$

\par In the following we assume for simplicity that $n=2m$ is even.  Let now
$a=\exp(\sum_{j=1}^m z_j e_j)$ with $z_j\in \R + i]-{\pi\over 4}, {\pi\over
4}[$.  Then by a straightforward calculation one gets that $$a(0)=J_{\tanh
z_1, \ldots, \tanh z_m}$$ and so $$T_\Omega(0)=\{ J_{z_1, \ldots, z_m}\:
z_j\in \C, \ |z_j|<1\}.$$ As a maximal abelian subspace $\e\subeq \q\cap
\p$ we choose $$\e=\{ \pmatrix{ 0 & iJ_{t_1,\ldots, t_m}\cr iJ_{t_1,\ldots,
t_m} & 0\cr}\:  t_1,\ldots,t_m\in \R\}.$$ Then one easily shows that
$$\exp(\e)(0)=\{ i J_{x_1, \ldots, x_m}\:  x_j\in \R, \ |x_j|<1\}$$ and in
particular $\exp(\e)(0)\subeq A\exp(i\Omega_0)(0)$.

\def\entries{

\[AG90 Akhiezer, D.\ N., and S.\ G.\ Gindikin, {\it On Stein extensions of
real symmetric spaces}, Math.\ Ann.\ {\bf 286}, 1--12, 1990

\[B03 Barchini,L., {\it Stein extensions of real symmetric spaces and the 
geometry of the flag manifold}, Math. Ann. {\bf 326} (2003), no. {\bf 2}, 331--346

\[BHH03 Burns, D., S.\ Halverscheid, and R.\ Hind, {\it The Geometry of
Grauert Tubes and Complexification of Symmetric Spaces}, Duke Math. J. {\bf 118}
(2003), no. {\bf 3}, 465--491  

\[F02 Faraut, J., {\it Formule de Gutzmer pour la complexification 
d'un espace riemannien sym\'etrique}, Rend. Mat. Acc. Lincei {\bf 9} (2002), 
vol. {\bf 13}, 233--241 

\[F03 Faraut, J., {\it Analysis on the crown of a Riemannian 
symmetric space}, in: Lie groups and symmetric spaces, 
Amer. Math. Soc. Transl. Ser. {\bf 2}, {\bf 210}

\[FK94 Faraut, J., and A.  Koranyi, ``Analysis on symmetric cones", Oxford
Mathematical Monographs, Oxford University Press, 1994

\[FT99 Faraut, J., and E.  G.  F.  Thomas, {\it Invariant Hilbert spaces of
holomorphic functions}, J.  Lie Theory {\bf 9} (1999), no.  {\bf 2}, 383--402

\[G68 Gangolli, R., {\it Asymptotic behaviour of spectra of compact
quotients of certain symmetric spaces}, Acta Math.  {\bf 121} (1968),
203--228

\[GK02a Gindikin, S., and B.\ Kr\"otz, {\it Complex crowns of Riemannian
symmetric spaces and non-compactly causal symmetric spaces}, Trans. Amer.
Math. Soc. {\bf 354 (8)} (2002), 3299--3327

\[GK02b ---, {\it Invariant Stein domains in Stein symmetric spaces and a
non-linear complex convexity theorem}, IMRN {\bf 18} (2002), 
959--971 

\[GK\'O03 Gindikin, S., B.\ Kr\"otz, and G.\ \'Olafsson, {\it Hardy spaces for
non-compactly causal symmetric spaces and the most continuous spectrum}, Math. Ann.
{\bf 327} (2003), 25--66

\[GK\'O04 ---,  {\it Holomorphic $H$-spherical distribution vectors 
in principal series representations}, Invent. math. {\bf 158} (2004), 
no. 3,   643--684

\[GM03 Gindikin, S., and T. Matsuki, {\it Stein Extensions of Riemann Symmetric 
Spaces and Dualities of Orbits in Flag Manifolds}, Transformation Groups {\bf 8} (2003), 
no. {\bf 4}, 333--376

\[He94 Helgason, S., ``Geometric Analysis on Symmetric Spaces'', 
Math.  Surveys and Monogr.  {\bf 39}, AMS, 1994

\[H\'O96 Hilgert, J.\ and G.\ \'Olafsson, ``Causal Symmetric Spaces,
Geometry and Harmonic Analysis,'' Acad.  Press, 1996

\[H\"o73 H\"ormander, L., ``An introduction to complex analysis in several
variables", North-Holland, 1973

\[Hu02 Huckleberry, A., {\it On certain domains in cycle spaces of flag 
manifolds}, Math. Ann. {\bf 323} (2002), no. {\bf 4}, 797--810 

\[K-T78 Kashiwara, M., A.\ Kowata, K.\ Minemura, K.\ Okamoto, T.\ Oshima,
and M.\ Tanaka, {\it Eigenfunctions of invariant differential operators on a
symmetric space}, Ann.  of Math.  {\bf (2) 107} (1978), no.  {\bf 1}, 1--39

\[K99 Kr\"otz, B., {\it The Plancherel theorem for biinvariant Hilbert
spaces}, Publ.  Res.  Inst.  Math.  Sci.  {\bf 35} (1999), no.  {\bf 1},
91--122

\[K01 --- , {\it Formal Dimension for Semisimple Symmetric Spaces}, Comp.
Math.  {\bf 125} (2001), 155-19

\[KN\'O97 Kr\"otz, B., K.--H.\ Neeb, and G.\ \'Olafsson, {\it Spherical
Representations and mixed symmetric spaces}, Representation Theory {\bf 1}
(1997), 424--461

\[KS04 Kr\"otz, B., and R.J.  Stanton, {\it Holomorphic extensions of
representations:  (I) automorphic functions}, Ann. Math. {\bf 159} (2004), 
no. {\bf 2}, 1--84

\[M03 Matsuki, T., {\it Stein extensions of Riemannian symmetric 
spaces and its generalization}, J. Lie Theory {\bf 13} (2003), no. {\bf 2}, 
565--572 

\[N65 Nagano, T., {\it Transformation groups on compact symmetric spaces}, 
Trans. AMS {\bf 118} (1965), 428-453

\[\'O91 \'Olafsson, G., {\it Symmetric Spaces of Hermitian Type},
Differential Geometry and its Applications {\bf 1}(1991), 195--233

\[S84 Schlichtkrull, H., ``Hyperfunctions and Harmonic Analysis on
Symmetric \linebreak Spaces", Progress in Math.  84, Birkh\"auser, 1984

\[St99 Stenzel, M., {\it The Segal-Bargmann Transform on a Symmetric Space 
of Compact Type}, J. Funct. Anal. {\bf 165} (1999), 44--58

} {\sectionheadline{\bf References} \frenchspacing \entries\par} \lastpage
\bye \end